\documentclass[11pt, nonatbib]{elsarticle}

\usepackage{appendix}
\usepackage[utf8]{inputenc}
\usepackage[T1]{fontenc}
\usepackage{amsmath,amsfonts,amssymb,amsthm}
\usepackage{algorithm}
\usepackage[noend]{algpseudocode}
\usepackage{graphicx}
\usepackage{caption}
\usepackage{subcaption}
\usepackage{natbib}
\usepackage{diagbox}
\usepackage{url}

\newtheorem{proposition}{Proposition}

\begin{document}

\title{Dynamic unsplittable flows with path-change penalties: new formulations and solution schemes for large instances}

\author[1]{François Lamothe\corref{cor1}}
\ead{francois.lamothe@isae-supaero.fr}
\author[1]{Emmanuel Rachelson}
\author[1]{Alain Haït}
\author[2]{Cédric Baudoin}
\author[3]{Jean-Baptiste Dupé}

\cortext[cor1]{Corresponding author}
\address[1]{ISAE-SUPAERO, Université de Toulouse, France}
\address[2]{Thales Alenia Space, Toulouse, France}
\address[3]{Centre national d'études spatiales (CNES), Toulouse, France}

\begin{abstract}
In this work, we consider the dynamic unsplittable flow problem. This variation of the unsplittable flow problem has received little attention so far. The unsplittable flow problem is an NP-hard extension of the multi-commodity flow problem where each commodity sends its flow on only one path. In its dynamic version, this problem features several time steps and a penalty is paid when a commodity changes its path from one time step to the next. We present several mixed-integer linear programming formulations for this problem and compare the strength of their linear relaxation. These formulations are embedded in several solvers which are extensively compared on small to large instances. One of these formulations must be solved through a column generation process whose pricing problem is more difficult than those used in classical flow problems. We present limitations of the pricing schemes proposed in earlier works and describe two new schemes with a better worst-case complexity. Overall, this work lays a strong algorithmic baseline for the resolution of the dynamic unsplittable flow problem, proposes original formulations, and discusses the compared advantages of each, thus hopefully contributing a step towards a better understanding of this problem for both OR researchers and practical applications.
\end{abstract}


\maketitle

\section*{Acknowledgments}
This document is the result of a research project funded by the Centre National d'\'Etudes Spatiales (CNES) and Thales Alenia Space.

\section{Introduction}

The unsplittable flow problem is a well-known variant of the multi-commodity flow problem. In this problem, one is given a directed or undirected graph together with capacities on its arcs, and a family of commodities, each characterized by an origin, a destination, and a demand. Each commodity has to route its demand from its origin to its destination through a unique path. The routing must ensure that capacities on the arcs are not exceeded by the flow of the commodities (or at least minimize the capacity overflow).

In this work, we take interest in a dynamic version of this problem which appears in several applications. Our practical application is the routing of communications inside constellations of satellites which is a central concern for the increasingly complex telecommunication constellations currently under development. In this context, a constellation manager would like to route the internet throughput of its users through the constellation while respecting the capacities of the telecommunication links. Because the constellation is orbiting the earth, its topology may be slightly changing over time and the users on the ground will not always be able to connect to the same satellites. Additionally, any modification to a user's path results in traffic disruption. Such path changes must thus be minimized. This dynamic problem can be modeled by a sequence of unsplittable flow problems representing different time steps. Each time step introduces a few changes to the underlying graph; some arcs are added or deleted hence changing the possible paths through the graph. Moreover, the commodities also become dynamic as their characteristics (origin, destination, and demand) may change over time. However, each commodity's characteristics and its time-dependent evolution remain associated with the same real-world object, e.g. the same user in the context of satellite communications. This sequence of unsplittable flow problems features two possibly conflicting objectives: 1) minimize the flow exceeding the capacities (overflow), 2) the commodities have to use the same path for consecutive time steps if possible (i.e. minimize the number of path changes over time). In applications, this problem might feature a very large number of time steps. The problem may then be approached through a rolling horizon by considering only a few time steps at a time. Because the previous horizons fixed a path for each commodity for the time-step preceding the current horizon and because these paths must be kept if possible to avoid path changes, each instance of the dynamic unsplittable flow problem is given with an initial path for each commodity. This initial path corresponds to the path used in the last time-step preceding the current horizon.

Our contributions are summarized as follows:
\begin{itemize}
    \item In the literature, \citet{gamvros2012multi} solved the dynamic unsplittable flow problem using a Mixed Integer Linear Program (MILP) inside a Branch and Price approach. After recalling the principle of their approach, we highlight that their pricing scheme has an exponential worst-case complexity. We then describe two new resolution approaches to the pricing problem with a polynomial worst-case complexity.
    \item We introduce two new formulations for the dynamic unsplittable flow problem. The first one is a mixed integer linear program that is shown to have a linear relaxation as strong as the formulation of \citet{gamvros2012multi}. This formulation yields good empirical results when the paths allowed for each commodity are heuristically restrained to a small carefully chosen set. The second formulation is a linear program that uses and links two types of flow variables that are not commonly seen at the same time in flow formulations: arc variables satisfying flow conservation constraints and path variables. This formulation can be rather competitive because the arc variables can be aggregated for commodities that share the same origin. However, because of this aggregation, this formulation can only be used to solve the linear relaxation of the dynamic unsplittable flow problem when only one time step is considered at a time (one decision time step where paths used in the previous time step are favored). Nevertheless, this formulation is used with a heuristic to obtain integer solutions of the problem, and yields, in our tests, better results than the formulations relying only on path variables.
    \item Finally, we describe several solvers for the dynamic unsplittable flow problem that use the presented formulations directly in MILP solvers or as a basis for heuristics. These solvers are extensively compared on a test-bed of small to large instances and several key aspects of the solvers explaining their performance are highlighted.
\end{itemize}

This paper is structured as follows. Related works on unsplittable flow problems are presented in Section \ref{sec:related_work}. In Section \ref{sec:dynamic_UFP} we recall the formulation of \citet{gamvros2012multi} and describe the two new formulations. One of them is shown to have a linear relaxation as strong as the one of \citet{gamvros2012multi}. In Section \ref{sec:solving_path-sequence_formulation} we explain how the formulation of \citet{gamvros2012multi} may be solved through column generation. Then, after recalling their pricing scheme for this formulation, we highlight its shortcomings and describe two new pricing schemes that have a lower worst-case complexity. Section~\ref{sec:solving_dynamic_UFP} discusses how the formulations presented in the previous sections can be embedded in solvers. Computational results for these solvers are presented in Section~\ref{sec:computational_results}. Finally, we conclude in Section~\ref{sec:conclusion}.

\section{Related work}
\label{sec:related_work}

To the best of our knowledge, the dynamic version of the unsplittable flow problem with path change penalties has only been studied in the work of \citet{gamvros2012multi}. They proposed a Branch and Price and Cut algorithm based on a new formulation named the path-sequence formulation. Their main contribution, besides this new formulation, is a pricing problem, which enables them to solve to optimality the linear relaxation of the path-sequence formulation. Although they studied a problem also containing network design features (decisions need to be taken about which arcs or arc configurations are allowed), their algorithm is valid for the dynamic unsplittable flow problem. They showed that their approach provides up to $13\%$ better results than one-step lookahead methods when exact solvers are used in both cases.

Other similar problems have been studied in their dynamic variant. For instance, approaches for solving multiple time steps network design have been proposed by \citet{lee2009dynamic, contreras2011dynamic, fragkos2017multi}. In these approaches, opening and closing of arcs and nodes of the network has a cost. These costs are very similar to the path change penalties studied in the dynamic unsplittable flow problem. Similarly, vehicle routing problems have a variant where serving a user with the same vehicle over time is valued \citep{kovacs2014vehicle, luo2015service, stavropoulou2019vehicle}. This persistence of vehicle choice over time recalls the incentive to keep the same path in the dynamic unsplittable flow problem.

Contrarily to its dynamic variant, the static unsplittable flow problem has been largely studied. This section does not aim at providing an extensive overview of this literature but focuses on recalling the most successful approaches among exact methods, approximation algorithms, and meta-heuristics. The seminal exact method is the one of \citet{barnhart2000using}. They presented a Branch and Price and Cut procedure applied to a formulation called the path formulation. Most posterior works use this baseline as a comparison. A major contribution of their work is their branching strategy which is used in most subsequent exact methods such as those of \citet{alvelos2003comparing, park2003integer, gamvros2012multi}. They also included lifted cover inequalities of the capacity constraints to strengthen the linear relaxation of their formulation. \citet{park2003integer} mixed a path formulation and a knapsack formulation to derive a new linear formulation of the problem. The linear relaxation of this formulation yields a stronger lower bound which in turn decreases the time needed to complete the branching procedure. \citet{belaidouni2007minimum} presented a cutting plane method based on superadditive functions to generate strong cuts for their Branch and Price method. It appears in small instances that the addition of such cuts provides integer solutions without using a Branch and Bound procedure. Finally, \citet{fortz2016convex, fortz2017models} studied an unsplittable flow problem where the cost of using each arc is convex and piece-wise linear which is also the case for the costs used in our formulation.

The approximation properties of the static unsplittable flow problem are also well-known. We recall here some results when the objective function is the congestion, \textit{i.e.} the smallest number by which it is necessary to multiply all the capacities to fit all the commodities. This objective function is the closest to the one studied in this work. We refer to the Handbook of approximation algorithms \citep{gonzalez2007handbook} for a detailed survey on approximation algorithms in the context of unsplittable flows. The algorithm with the best approximation factor for congestion is a randomized rounding method introduced by \citet{raghavan1987randomized}. This procedure is a $O\left(\frac{\ln |E|}{\ln \ln |E|}\right)$-approximation algorithm that works for directed and undirected graphs, where $|E|$ is the number of arcs in the network. This algorithm was extended by \citet{rivano2002lightpath} even though the approximation properties were not proven for the extended algorithm. \citet{lamothe2021randomized} further extended the previous randomized rounding algorithm to yield better results in practice while keeping the same approximation factor. This algorithm is called Sequential Randomized Rounding (SRR) and will be used in this work to create integer solutions from the linear relaxation of randomized rounding algorithms.

Finally, various meta-heuristics were investigated in the unsplittable flow literature. We refer to the work of \citet{li2010ant} and \citet{santos2013hybrid} for the most efficient meta-heuristics. \citet{li2010ant} showed that their ant colony optimization method compared favorably with the CPLEX solver, and was able to solve instances with up to 60 nodes, 400 arcs, and 3500 commodities to optimality in less than 900 seconds. \citet{santos2013hybrid} created an algorithm that solves the linear relaxation of the problem with column generation and uses the paths created during the column generation inside a Greedy Randomized Adaptive Search Procedure (GRASP). They solve all their instances (26 nodes, 80 arcs, 500 commodities) in less than 180 seconds with values close to the linear relaxation's lower bound.

\section{Formulations for the dynamic unsplittable flow problem}
\label{sec:dynamic_UFP}

The dynamic unsplittable flow problem is a variation of the unsplittable flow problem where several time steps are considered. This variation features two conflicting objectives: respecting as much as possible the arc capacities and minimizing the number of times each commodity changes the path it uses. The path changes are modeled with penalties. For the sake of simplicity, we consider that penalties are uniform over all commodities and time steps. The price of a penalty will be denoted $\alpha$ and can be used as a scaling factor between the two objective functions. The extension to the general case of heterogeneous penalties does not involve major changes in the algorithms and formulations. 

Note that the objective function used in this work differs slightly from the objective function classically used for the unsplittable flow problem. Herein, we focus on minimizing the violation of the capacities. In this context, one classically minimizes the congestion which is the smallest number by which it is necessary to multiply all the capacities in order to fit all the commodities. Minimizing the congestion puts no restrictions on the flow going through the arcs that do not have maximal congestion. In particular, it induces no incentive to minimize the congestion on those arcs. This becomes problematic when an arc is largely more congested in every solution than any other arc because it lifts all restrictions for the other arcs. Since in our constellation application, we would like to minimize the violation on all parts of the graph independently of the congestion of the other parts, we instead consider an objective function penalizing the sum of overflow on the arcs. Moreover, it was shown by \citet{lamothe2021randomized} that this objective function enables the SRR heuristic to create solutions with both a lower cumulated overflow and better congestion.
We consider that, for each time step, the cumulated overflow over the arcs is not penalized until it reaches a threshold denoted $B$. After this allowed amount, the overflow is heavily penalized compared to the path changes (because the penalty price $\alpha$ of a path change is set to a small value). This encourages the algorithms to find solutions with a minimum number of path changes among the ones with less than $B$ cumulated overflow. Note that no preference is given on how the overflow is distributed among the arcs as long as its sum in each time step is lower than the threshold $B$. In general, any piecewise linear penalization of the overflows can be considered with the methods proposed in this work, and extension to the general case is quite straightforward.

In this section, we present a variation of the formulation introduced by \citet{gamvros2012multi} to our specific version of the dynamic unsplittable flow problem and present a new MILP formulation that is shown to have a linear relaxation as strong as theirs. Then, we present a new linear programming formulation that models the linear relaxation of the dynamic unsplittable flow problem when only one time step is considered. This formulation uses an aggregation of the commodities that drastically reduces its number of variables. However, because of this aggregation, it cannot be transformed into a mixed-integer linear formulation of the dynamic unsplittable flow problem. It will thus be used in our experiments inside a heuristic that creates integer solutions from the linear relaxation of unsplittable flow problems. For the remainder of this paper, we will use the following notations where $t \in T$ is the time step index and $k \in K$ is the index of commodities:
\begin{itemize}
    \item $G = (V, E, (E_t)_{t \in T})$ is a directed or undirected dynamic graph where $V$ and $E$ are the sets of nodes and arcs and where $(E_t)_{t \in T}$ are the sets of allowed arcs at each time step;
    \item $(O^k_t, D^k_t, d^k)_{t \in T, k\in K}$ is a set of commodities defined by their origin node $O^k_t \in V$, destination node $D^k_t \in V$ at each time-step, and a fixed demand $d^k \in \mathbb{R}^+$;
    \item $(c_{et})_{e \in E_t, ~ t \in T}$ are capacities on the arcs of the dynamic graph at each time step.
\end{itemize}

\subsection{The path-sequence formulation}
\label{sec:path-sequence-formulation}

The following formulation is a variation of the one introduced by \citet{gamvros2012multi} to our specific version of the dynamic unsplittable flow problem. In the model of \citet{gamvros2012multi}, arc capacities cannot be exceeded while the use of each arc by a commodity has a cost. The objective is then to minimize the total induced cost. In our case, we allow the arc capacities to be exceeded and minimize some of the excess flow. The number of variables in this formulation is not polynomial, neither in the number of nodes/arcs nor in the number of time steps, and its linear relaxation must be solved through column generation.

In this formulation, the meaning of the variables is the following:
\begin{itemize}
    \item $x_s^k$ decides if  commodity $k$ uses the path-sequence $s$ over the considered horizon. If a path-sequence $s = (p_1, ..., p_{|T|})$ is used, then path $p_j$ is used at time step $j$. Each path-sequence induces a number of path-changes $n_s^k$;
    \item $o_{et}$ represents the overflow on arc $e$ at time step $t$;
    \item $o_t$ represents the amount of overflow that exceeds at time step $t$ the amount $B$ of unpenalized overflow.
\end{itemize}

We denote $S^k$ the set of all path-sequences usable by commodity $k$ while $S^k_{et}$ denotes the set of all path-sequences usable by commodity $k$ for which the path at time step $t$ is going through arc $e \in E_t$. The path-sequence formulation is given by:

\begin{subequations}
\begin{alignat}{2}
&\min_{x_s^k, o_{et}, o_t} ~ \alpha \sum_{k \in K} \sum_{s \in S^k} n_s^k ~ x_s^k + \sum_{t \in T} o_t & \qquad &  \label{eq: obj_e}\\
&\text{subject to} &      & \notag \\
&\sum_{s \in S^k} x_s^k = 1 & & \forall k \in K \label{eq: convex_e}\\
&\sum_{k \in K} \sum_{s \in S^k_{et}} x_s^k ~ d^k \leq c_e + o_{et} & & \forall e \in E_t, ~ \forall t \in T \label{eq: capacity_const_e}\\
&\sum_{e \in E_t} o_{et} \leq B + o_t & &  \forall t \in T \label{eq: capacity_const_tot_e}\\
& x_s^k \in \{0,1\} & & \forall s \in S^k, ~ \forall k \in K \label{eq: integrity_e} \\
& o_{et} \in \mathbb{R}^+, ~ o_t \in \mathbb{R}^+ & & \forall e \in E_t, ~ \forall t \in T
\end{alignat}
\end{subequations}

The objective function \eqref{eq: obj_e} is composed of two terms: the sum of the path-change penalties and the sum over the time steps of the overflow exceeding the threshold $B$. Note that even though there is a very large number of path-sequences, because each commodity is limited to only one path-sequence by the equations (\ref{eq: convex_e}) and each path-sequence induces at most $|T|-1$ path-changes, the number of path-change penalties is upper-bounded by $|K|(|T|-1)$. The equations (\ref{eq: capacity_const_e}) and (\ref{eq: capacity_const_tot_e}) are the capacity constraints; they ensure that $o_{et}$ represents the overflow on arc $e$ at time step $t$ and that $o_t$ represents the amount of overflow that exceeds the threshold $B$ at time step $t$. The variables $ x_s^k$ being binary in the equations \eqref{eq: integrity_e} ensures that the flow is unsplittable.

\subsection{Extended arc-path formulation}
\label{dynamic_path_formualtion}

The static unsplittable flow problem has two classical formulations, the arc-node formulation and the arc-path formulation, which can be found in \citep{barnhart2000using}. These two formulations can be extended to the dynamic unsplittable flow problem. For the arc-path formulation, the obtained formulation is presented below. This formulation uses a polynomial number of variables and constraints in the number of commodities and time steps but not in the number of nodes and arcs because of the exponential number of paths in a graph. The variables of this formulation have the following meaning:
\begin{itemize}
    \item $x_{pt}^k$ decides if commodity $k$ uses path $p$ at time step $t$;
    \item $o_{et}$ represents the overflow on arc $e$ at time step $t$;
    \item $o_t$ represents the amount of overflow that exceeds, at time step $t$, the amount $B$ of unpenalized overflow;
    \item $n_{pt}^k$ takes value one if commodity $k$ uses path $p$ at time step $t$ but not at time step $t-1$.
\end{itemize}

We denote $P_t^k$ the set of all the paths usable by commodity $k$ at time step $t$. 
The set of time steps $T$ is extended to include a time step 0 to take into account imposed initial commodity paths. The set $P_0^k$ of usable paths at time step zero contains only one path which is the one taken by commodity $k$ before the decision time interval. 
Note that the variables $x_{pt}^k$ are considered to exist for both admissible ($p\in P^k_t$) and inadmissible ($p \not\in P^k_t$) paths but a variable $x_{pt}^k$ is set to zero when path $p$ is inadmissible. The extended arc-path formulation is given by:

\begin{subequations}
\begin{alignat}{2}
&\min_{x_{pt}^k,n_{pt}^k, o_{et}, o_t} ~ \alpha \sum_{k \in K, t \in T} n_{pt}^k + \sum_{t \in T} o_t & \qquad &  \label{eq: obj}\\
&\text{subject to} &      & \notag \\
&\sum_{p \in P_t^k} x_{pt}^k = 1 & & \forall k \in K, ~ \forall t \in T \label{eq: convex}\\
&\sum_{k \in K} \sum_{p \in P_t^k |e \in p} x_{pt}^k ~ d^k \leq c_{et} + o_{et} & & \forall e \in E_t, ~ \forall t \in T \label{eq: capacity_const}\\
&\sum_{e \in E_t} o_{et} \leq B + o_t & &  \forall t \in T \label{eq: capacity_const_tot}\\
& x_{pt}^k - x_{p,t-1}^k \leq n_{pt}^k & & \forall p \in P_t^k, ~ \forall k \in K, ~ \forall t \in T \setminus \{0\} \label{eq: delta_1}\\
& x_{pt}^k \in \{0,1\}, ~ n_{pt}^k \in \mathbb{R}^+  & & \forall p \in P_t^k, ~ \forall k \in K, ~ \forall t \in T \label{eq: integrity}\\
& o_{et} \in \mathbb{R}^+, ~ o_t \in \mathbb{R}^+ & & \forall e \in E_t, ~ \forall t \in T
\end{alignat}
\end{subequations}
The objective function \eqref{eq: obj} is composed of two terms: the sum of the path-change penalties and the sum over the time steps of the overflow exceeding the threshold $B$. Note that even though there is a large number of paths at each time-step, because each commodity is limited to only one path per time-step by the equations (\ref{eq: convex}), the number of path-change penalties is upper-bounded by $|K|(|T|-1)$. The equations (\ref{eq: capacity_const}) and (\ref{eq: capacity_const_tot}) are the capacity constraints: they ensure that $o_{et}$ represents the overflow on arc $e$ at time step $t$ and that $o_t$ represents the amount of overflow that exceeds the threshold $B$ at time step $t$. The equations (\ref{eq: delta_1}) ensure that $n_{pt}^k$ takes value 1 when a change of path occurs for commodity $k$ between time step $t$ and time step $t-1$. The variables $x_{pt}^k$ being binary in the equations \eqref{eq: integrity} ensures that the flow is unsplittable.

The above formulation is an adaptation of the arc-path formulation for unsplittable flows to the multi-timestep setting with path change penalties. It is to be noted that a similar adaptation can be achieved with an arc-node formulation. The resulting formulation is compact which enables an exact resolution method with a commercial solver. However, due to the very high number of variables of this compact formulation, this method does not scale well even in relatively small instances. The extended arc-node formulation is given in \ref{appendix_arc_node} together with the computational results for small instances.

\subsection{Equivalence of the relaxation of the two MILP models}

Combinatorial problems can often be mathematically described using several MILP models. These models are usually compared through their number of variables and number of constraints but also through the strength of their linear relaxation. Indeed, the stronger the linear relaxation, the more suited the MILP model is for a resolution using a Branch-and-Bound procedure, and the more information the linear relaxation provides on the integer problem. In this section, we show that the two MILP formulations presented above for the dynamic unsplittable flow problem have equally strong linear relaxations.

\begin{proposition}
Let $V_1$ be the value of the linear relaxation of the extended arc-path formulation and $V_2$ be the value of the linear relaxation of the path-sequence formulation, then $V_1 = V_2$.
\end{proposition}

Proof of the above proposition is given in \ref{appendix_proof}. The two formulations can thus be used interchangeably and the main criterion for choosing one over the other should be its resolution time.

\subsection{Aggregated arc-node formulation for one time step}
\label{sec:arc_node_one_time step}

The aggregated arc-node formulation below models the linear relaxation of a dynamic unsplittable flow problem (which is the same problem except that each commodity can use several paths to route its flow since the integrality of the variables is relaxed) when one time step is considered. It uses and links two types of flow variables that are not commonly seen at the same time in flow formulations: arc variables satisfying flow conservation constraints that are common in arc-node formulations and path variables that are used in arc-path formulations. This formulation is especially relevant because path-change penalties are better represented with path variables while flows can be efficiently represented with arc-node formulations when commodities are grouped by origin to create super-commodities. A super-commodity $k'$ contains several commodities $k$, and its demand is the sum of their demands between their common origin and their destinations. Because this formulation considers splittable flows and because the flow of each commodity is indistinguishable except for their origin and destination, grouping the commodities can be done without loss of generality and without changing the flow distribution represented by the optimal solution. However, it forces the formulation to provide aggregated flow. Thus, after solving the formulation, it is necessary to compute which exact paths are used by each commodity. This can be done exactly and quickly in $O(|V|(|E| + |K|))$ operations with a flow decomposition algorithm \citep{ford1956network}. Moreover, grouping the commodities greatly reduces the number of variables of the arc-node formulations and enables the computation of linear relaxations for large instances that would otherwise be intractable with such formulations. 

Note however that, because of this grouping, the linear formulation below cannot be converted into a mixed-integer linear program representing the dynamic unsplittable flow problem. Indeed, the grouping prevents the formulation from choosing a unique path separately for each commodity which is necessary to ensure that the flow of each commodity is unsplittable. Moreover, considering only one time step in the formulation (one decision time step where paths used in the previous time step are favored) enables the use of only one path variable per commodity; a simplification that cannot be generalized to longer time horizons. This simplification allows the formulation to have a polynomial number of variables in the number of nodes and arcs of the graph. In compensation, the formulation can only model the linear relaxation of a dynamic unsplittable flow problem when one time step is considered. In order to create unsplittable solutions to the dynamic unsplittable flow problem, this formulation will be used in conjunction with the SRR heuristic of \citet{lamothe2021randomized}. For the sake of readability, we will drop for this section the index $t$ indicating the considered time step. The variables of this formulation have the following meaning:
\begin{itemize}
    \item $f_e^{k'}$ indicates how much flow the super-commodity $k'$ pushes on arc $e$;
    \item $x_{p^k}^k$ decides the proportion of the flow of commodity $k$ which is sent on the path $p^k$ used in the previous time step;
    \item $o_e$ represents the overflow on arc $e$;
    \item $o$ represents the amount of overflow that exceeds the amount $B$ of unpenalized overflow.
\end{itemize}

Since only one time step is considered, if $p^k$ is the path used in the previous time step by commodity $k$, in the current time step the number of path changes is equal to $1 - x_{p^k}^k$. Thus, the path changes can be taken into account in the objective by the term $\sum_{k \in K} (1 - x_{p^k}^k)$. In this formulation we also use the Kronecker notation: $\delta^x_y = 1$ if $x = y$ and $\delta^x_y = 0$ otherwise. Moreover, $E^-(v)$ and $E^+(v)$ will denote the set of incoming and outgoing arcs of a node $v$. The aggregated arc-node formulation is given by:

\begin{subequations}
\begin{alignat}{2}
&\min_{x_{p^k}^k, f_e^{k'}, o_e, o} ~ \alpha \sum_{k \in K} (1 - x_{p^k}^k) + o & \qquad &  \label{eq: obj_na_1}\\
&\text{subject to} &      & \notag \\
&\sum_{e \in E^+(v)} \hspace{-0.3cm} f_e^{k'} ~ \hspace{-0.2cm} - \hspace{-0.3cm} \sum_{e \in E^-(v)} \hspace{-0.3cm} f_e^{k'} = \sum_{k \in k'} d^k \delta^{O^k}_v - \sum_{k \in k'} d^k \delta^{D^k}_v & & \forall k' \in K', ~ \forall v \in V \label{eq: flow_const_na_1}\\
&\sum_{k \in k'| e \in p^k} x_{p^k}^k d^k \leq f_e^{k'} & & \forall k' \in K', ~ \forall e \in E \label{eq: linking}\\
&\sum_{k' \in K'} f_e^{k'} \leq c_e + o_e & & \forall e \in E \label{eq: capacity_const_na_1}\\
&\sum_{e \in E} o_e \leq B + o & & \label{eq: capacity_const_tot_na_1}\\
& f_e^{k'} \in \mathbb{R}^+, x_{p^k}^k \in [0,1], ~ o_e \in \mathbb{R}^+, ~ o \in \mathbb{R}^+ & & \forall k' \in K', ~ \forall e \in E, ~ \forall k \in K
\end{alignat}
\end{subequations}

The objective function \eqref{eq: obj_na_1} is composed of two terms: the sum of the path-change penalties and the overflow exceeding the threshold $B$. The equations (\ref{eq: flow_const_na_1}) is the flow conservation constraint for the super-commodity $k'$. The equations \eqref{eq: linking} are a linking constraint between the flow variables $f_e^{k'}$ and the path variables $x_{p^k}^k$. The equations (\ref{eq: capacity_const_na_1}) and (\ref{eq: capacity_const_tot_na_1}) are the capacity constraints: they ensure that $o_e$ represents the overflow on arc $e$ and that $o$ represents the amount of overflow that exceed the threshold $B$.

\section{Pricing schemes for the path-sequence formulation}
\label{sec:solving_path-sequence_formulation}

The path-sequence formulation was introduced by \citet{gamvros2012multi} to solve the dynamic unsplittable flow problem. Because of its large number of variables $x_s^k$, a column generation process is required to solve its linear relaxation. This process relies on a so-called pricing scheme that must be tailored to the problem at hand. In this section, we start by briefly recalling the concept of column generation and present how the pricing problem particularizes to the context of dynamic flows. An algorithm solving this problem is given in the work of \citet{gamvros2012multi} which relies on the computation of k-shortest paths. We partially recall their method and highlight the fact that this pricing scheme has an exponential worst-case complexity. This is mainly because there might exist an exponential number of paths of the same length/cost in a graph. To alleviate this problem, we propose two new pricing schemes that do not rely on k-shortest path computation and have a polynomial worst-case complexity.

\subsection{Column generation for the path-sequence formulation}
\label{sec:CG and path-sequence}

Column generation is a technique applied when the number of variables is very large. In this case, applying a standard simplex procedure and computing the reduced cost of each variable is intractable. However, it can be tractable to compute the variable with the smallest reduced cost through an optimization problem called the pricing problem. Recall that the reduced cost of a variable $x_i$ is $r_i = c_i - u \cdot A_i$ where $c_i$ is the cost of $x_i$ in the objective function, $u$ is the vector of dual variables (one dual variable per constraint) and $A_i$ is the vector of coefficients associated with $x_i$ in the constraints.

In the context of the path-sequence formulation, the cost of a variable $x_s^k$ in the objective function \eqref{eq: obj_e} is $n_s^k$ the number of path-changes induced by the path-sequence $s$. As for the constraints, considering only the constraints with non a zero coefficient, a variable $x_s^k$ has a unit coefficient in the constraint (\ref{eq: convex_e}) associated with commodity $k$. Moreover, it appears in several constraints of type (\ref{eq: capacity_const_e}). If the path-sequence $s$ uses an arc $e$ at time-step $t$ then $x_s^k$ appears in the associated constraint with coefficient $d^k$.

This means that a variable $x_s^k$ with $s = (p_1, ..., p_{|T|})$ has the following reduced cost: 
$$r_s^k = n_s^k - u^k - d^k ~ \sum_{t \in T} \sum_{e \in p_t} u_{et}$$
where $u^k$ is the dual variable of the constraint (\ref{eq: convex_e}) associated with commodity $k$, and $u_{et}$ is the dual variable of the constraint (\ref{eq: capacity_const_e}) associated with arc $e \in E_t$ at time step $t$. Note that the coefficients for the constraints \eqref{eq: convex_e} are the same for all path-sequences of the same commodity. Because, the methods used to solve the pricing problem search for the best path-sequence for each commodity separately, the term $- u^k$ is constant in the pricing problem of a commodity. It can thus be removed from the pricing problem and added separately to the reduced cost after the computation.

Thus, the problem of finding the path-sequence of smallest reduced cost for a commodity $k$ is the following: choose the path-sequence with the smallest cost such that choosing a path $p_t$ at time step $t$ induces a cost $-d^k ~ \sum_{e \in p_t} u_{et}$ and not keeping the same path from $t$ to $t+1$ induces a cost of one.

\subsection{The pricing method of \citet{gamvros2012multi} and its limitations} 
\label{sec:pricing_gamvros}

As shown by \citet{gamvros2012multi}, the pricing problem of the path-sequence formulation can be decomposed into two parts. First, determine a set of candidate paths $\hat{P}_t^k$ for each time step. Second, using the paths in $\hat{P}_t^k$, construct a path-sequence of smallest reduced cost. This second part can be solved through dynamic programming. Indeed, choosing the best path for a time step $t$ depends on the path chosen at time step $t-1$ but not on all the previously chosen paths. For more details, see the work of \citet{gamvros2012multi}. 

The set $\hat{P}_t^k$ of usable paths cannot be the entire set $P_t^k$ of valid paths because the cardinality of $P_t^k$ is exponential in the size of the graph. Thus, \citet{gamvros2012multi} restrain $\hat{P}_t^k$ to be the result of a k-shortest path algorithm on the graph $G_t = (V, E_t)$, where the cost of each arc $e \in E_t$ is equal to $- d^k u_{et}$. Let us denote the returned paths $ p_1, ..., p_{\kappa} $ where $ p_1 $ is the shortest path and $ p_{\kappa} $ the longest path returned by the algorithm. The number of paths $ \kappa $ calculated by the k-shortest path algorithm is fixed such that it is the maximum integer satisfying $ c(p_{\kappa}) \leq c(p_1) + 2 \alpha $ where $ c(p) $ denotes the cost of a path $ p $ and the value $2 \alpha$ is the cost of two path-change penalties. In other words, it is guaranteed that all the paths having a cost difference with the shortest path of less than $2 \alpha$ have been calculated. 

A sequence of paths of minimum reduced cost can be calculated using only the paths satisfying this condition. Indeed, suppose that a sequence of paths of minimum reduced cost uses at time step $ t $ a path $ p $ not satisfying $ c(p) \leq c(p_1) + 2 \alpha $. We can then replace $ p $ by $ p_1 $ to create a sequence of paths of lower cost. Indeed, the exchange can only increase by two the number of path changes made in the sequence of paths, and this increase is compensated by a decrease in the cost of the path by an amount of at least $2 \alpha$. The new sequence of paths is therefore of minimum reduced cost and uses one of the paths returned by the k-shortest path algorithm.
\textit{Remark: if the set of paths allowed at each time-steps for each commodity is restricted to a specific set (as this will be the case for several algorithms of the experimental section), then this restricted set can be used as the set of candidate paths $\hat{P}_t^k$ in the pricing problem. This makes the pricing problem easier to solve.}

In their work, \citet{gamvros2012multi} proposed enhancements of the above method to create the sets $\hat{P}_t^k$. However, these methods use a k-shortest path algorithm with a stopping condition based on the cost difference between the shortest and the longest computed path. These methods all have the following drawback: there might be an exponential number of paths whose cost difference with the shortest path is less than a constant. This may happen in at least the two following situations. First, the constant is the penalty induced by two path changes. If the price of a path change is rather large compared to the (reduced) costs of the paths then all the paths of the graph must be computed. Secondly, if the graph on which the k-shortest path computations are made contains a lot of zero-cost arcs then there might be a lot of equally shortest paths that must all be computed. In the limit where all arcs have a zero cost, all the paths of the graph must be computed. This second case tends to happen in our version of the path-sequence formulation if only a few arcs are overloaded. However, this was not the case in the tests of \citet{gamvros2012multi} because, in their formulation, the commodities have to pay a positive price to use an arc. Thus in their cases, the cost of each arc of the k-shortest path graph is positive. Nevertheless, their pricing algorithm has an exponential worst-case complexity.

\subsection{New pricing schemes without k-shortest paths}
\label{sec:new_pricing}

In this section, we introduce two new pricing schemes that do not rely on k-shortest path computations and have a polynomial worst-case complexity.

\textbf{Pricing with only shortest paths.} In a path-sequence, let us denote $t_1 ... t_2$ any sequence of consecutive time steps such that the path-sequence changes its path before $t_1$ and after $t_2$, but not between $t_1$ and $t_2$. Let $p_{t_1 ... t_2}$ be the path yielding the lowest total reduced cost on this sequence of time steps if it exists. The method is based on the following statement. The path $p_{t_1 ... t_2}$ is independent of the paths chosen outside of $t_1 ... t_2$ and can be computed with a single call to a shortest path algorithm on a graph $G_{t_1 ... t_2}$. The nodes of $G_{t_1 ... t_2}$ are $V_{t_1 ... t_2} = V$, its arcs are $E_{t_1 ... t_2}  = \bigcap_{t = t_1}^{t_2} E_t$ and the cost on an arc $e \in E_{t_1 ... t_2}$ is $c_e = \sum_{t = t_1}^{t_2} u_{et}$. Note that there may not exist a valid path for a commodity on $G_{t_1 ... t_2}$. In that case, the commodity must change its path between $t_1$ and $t_2$. If $\hat{P}_t^k$ is set to $\{p_{t_1 ... t_2}| ~ \forall t_1 ... t_2 \text{ such that } t_1 \leq t \leq t_2\}$ then a path-sequence of smallest reduced cost can be computed using only the paths in $\hat{P}_t^k$. Indeed, suppose that a path-sequence of smallest reduced cost uses at time step $t$ a path $p$ that is not in $\hat{P}_t^k$. Path $p$ is kept on a sequence $t_1 ... t_2$ and can be replaced by $p_{t_1 ... t_2}$ on this sequence to create a path-sequence of smaller reduced cost.

We now discuss the complexity of the pricing method presented above. First, we consider the complexity of computing the paths $p_{t_1 ... t_2}$. Let us denote $|E| = \max_t |E_t|$. For a sub-sequence $t_1 ... t_2$, the method starts by computing $G_{t_1 ... t_2}$. Note that $G_{t_1 ... t_2}$ has the same nodes as $G_{t_1 ... t_2-1}$ and its arcs are $E_{t_1 ... t_2} = E_{t_1 ... t_2-1} \cap E_{t_2}$. Thus, except when $t_1 = t_2$ in which case $G_{t_1 ... t_2} = (V, E_{t_1})$, the graph $G_{t_1 ... t_2}$ can be computed from $G_{t_1 ... t_2-1}$ and $G_{t_2}$ in $O(|E_{t_1 ... t_2-1}|)$ time which is equal to $O(|E|)$ time. Using a Dijkstra algorithm, the computation of $p_{t_1 ... t_2}$ from $G_{t_1 ... t_2}$ takes $O(|E_{t_1 ... t_2}| ~ \ln|E_{t_1 ... t_2}|)$ time which is also $O(|E| ~ \ln|E|)$ time. Overall, as there are $\frac{|T| (|T| - 1)}{2}$ paths of this form to compute, the complexity of computing all the paths $p_{t_1 ... t_2}$ is $O(|T|^2 |E|\ln|E|)$. Compared to the method of \citet{gamvros2012multi}, the number of paths computed is quadratic in the number of time steps instead of linear. However, in this method, the number of paths computed is constant in the size of the graph, while the method of \citet{gamvros2012multi} might compute an exponential number of paths in the size of the graph. The second part of the pricing method is the dynamic programming algorithm of \citet{gamvros2012multi}. It has a complexity of $O(\sum_{t \in T} |\hat{P}_t^k|)$ which is equal in our case to $O(\sum_{t \in T} t(|T|-t))$ which is equal to $O(|T|^3)$. Thus the overall complexity of the above method is $O(\max\{|T|^3, |T|^2 |E|\ln|E|\})$. We now describe an all-in-one method to solve the pricing problem that does not use the dynamic programming algorithm of \citet{gamvros2012multi} and has an even better complexity.

\begin{figure}
     \centering
     \includegraphics[width=0.7\textwidth]{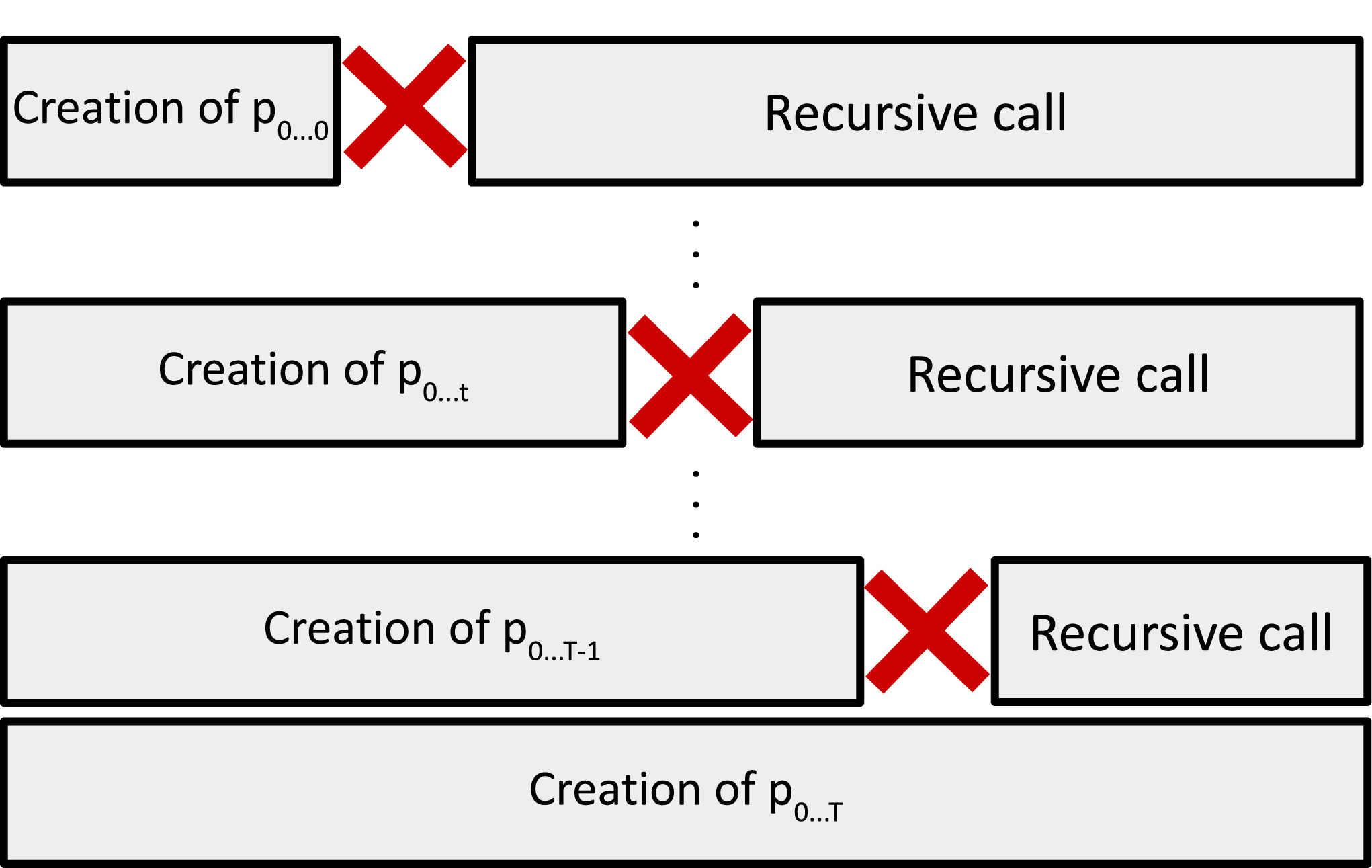}
     \caption{Illustration of the case disjunction made by the All-in-one pricing on the position of the first path change.}
     \label{case_disjunction}
\end{figure}

\textbf{Pricing all-in-one.} This method uses the paths $p_{t_1 ... t_2}$ as defined and constructed in the paragraph "Pricing with only shortest paths". This method makes a case disjunction on the position of the first path change. Indeed, if the first path-change occurs after just time step $t$ then the best path-sequence can be created by using the path $p_{0 ... t}$ for the time steps $0$ to $t$ and juxtaposing the optimal path-sequence for the sub-sequence of time steps $t+1, ..., |T|$. The latter optimal path-sequence can be computed by a recursive call of the method on the sub-sequence of time steps $t+1, ..., |T|$. Note that this recursive call might have already been made in another part of the algorithm. In this case, the algorithm does not need to re-make the computation and can just use the already computed path-sequence. Once all positions for the first time step have been considered, the optimal path-sequence for the whole horizon is the best among the ones created in each case. The case disjunction of this method is illustrated in Figure \ref{case_disjunction}.

We now study the complexity of this last method. First, from the complexity study of the previous pricing scheme, we know that all the paths $p_{t_1 ... t_2}$ can be computed in $O(|T|^2 |E| \ln|E|)$ time. Secondly, there are $|T|+1$ (recursive) calls to the algorithm, one for each sub-sequence of time steps of the form $t ... |T|$. Each call considers $|T| - t$ cases that are treated in constant time if we assume the paths $p_{t_1 ... t_2}$ and the recursive calls for smaller sub-sequences of time steps have already been computed. Thus all the recursive calls and case disjunctions are made in $O(|T|^2)$ time. The total complexity of the method is thus $O(|T|^2 |E| \ln|E|)$.

For the clarity of the presentation, we described all the pricing methods as if no initial path was given to each commodity. To account for the initial paths, just add the initial paths in the set of considered paths $\hat{P}_t^k$ in the methods relying on \citet{gamvros2012multi} dynamic programming algorithm, and in the last method make an additional case disjunction on how long the initial path should be kept.

\section{Solving the dynamic unsplittable flow problem}
\label{sec:solving_dynamic_UFP}

In the previous sections, we presented several formulations for the dynamic unsplittable flow problem and new ways to solve the linear relaxation of the path-sequence formulation. We now describe how to use these formulations to create resolution methods for the problem. In particular, most of the solvers presented in the experimental section create integer solutions of the problem from its linear relaxation using the Sequential Randomized Rounding heuristic (SRR) which is presented in the following sub-section. We then discuss the concept of restriction of the set of usable paths in a solver. Finally, we summarize the solvers that will be compared in the experimental section.

\subsection{The SRR heuristic}
\label{sec:SRR_heuristic}

Proposed and analyzed by \citet{lamothe2021randomized}, the Sequential Randomized Rounding heuristic (SRR) is a greedy heuristic shown to create very good solutions of the static unsplittable flow problem from its linear relaxation. Thus we will apply it to various formulations to create solutions of the dynamic unsplittable flow problem. The SRR heuristic initially computes the linear relaxation of the problem and then alternates between two actions to find an integer solution:
\begin{itemize}
    \item select the unfixed commodity with the largest demand; use randomized rounding to fix this commodity to a single path;
    \item compute the linear relaxation of the modified problem where the chosen commodities have been fixed to their single path;
\end{itemize}

The randomized rounding step was introduced by \citet{raghavan1987randomized}. It consists in choosing a path for a commodity according to the information given in the linear relaxation. In the linear relaxation, a proportion $x_p^k$ of the flow of commodity $k$ is sent on path $p$. A randomized rounding step consists in fixing commodity $k$ to path $p$ in the integer solution with probability $x_p^k$.

The linear relaxation is updated at a certain frequency which is controlled by a hyper-parameter $\theta$. In the solution of the linear relaxation, some commodities use multiple paths. After $\theta$ of these commodities are fixed to a single path, the linear relaxation is updated.

In the case of the path-sequence formulation, note that a blunt application of the heuristic would recommend making all the necessary column generation iterations to find the optimal linear solution before applying randomized rounding steps. However, as the SRR algorithm is a heuristic, it is reasonable to work with an approximate solution of the linear relaxation. This means that we can perform only a few column generation iterations before applying randomized rounding steps. In return, we can update the linear relaxation more often. At first, a large number of column generation iterations are made to obtain a good approximation of the linear relaxation. Then, the SRR heuristic is applied by alternating between two types of actions:
\begin{itemize}
    \item update the linear relaxation through column generation iterations;
    \item choose a path for a commodity and a time step through randomized rounding.
\end{itemize}

\subsection{Restriction of the usable paths}
\label{sec:restriction}

At each time step, a large number of paths is usable by each commodity. If all these paths are considered, the dynamic unsplittable flow problem is challenging to solve due to its huge solution space. However, it is not necessary to consider all the paths to obtain good-quality solutions. Restricting the allowed paths per commodity to a small predefined set simplifies the problem. While the optimal solution might be lost, it is easier to find optimal solutions to the restricted problem. The restricted sets of paths used in this work are computed as follows. Let $P^k_{t\kappa}$ be the set of the $\kappa$-shortest paths for commodity $k$ at time step $t$ when the length of the arcs is equal to one (\textit{i.e.} with the least number of hops). For methods considering a single time step, the restricted set of usable paths is set to $P^k_{t\kappa}$. For methods considering multiple time steps, the restricted set of usable paths for commodity $k$ at time step $t$ is set to $\bigcup_{t' \in T} P^k_{t'\kappa} \cap P_t^k$ where $P_t^k$ is the set of all valid paths for commodity $k$ at time step $t$. Moreover, when a time step or a sequence of time steps is considered, the initial paths are also allowed in each time step where they are valid. This restriction has already been used in meta-heuristics designed for the unsplittable flow problem by \citet{laguna1993bandwidth} and \citet{masri2015multi}. Moreover, it allows us to use Branch and Bound algorithms without resorting to Branch and Price. The impact of this restriction in terms of solution quality and computing time will be evaluated in Section \ref{sec:computational_results}.

\subsection{Solvers}

In the following, we present all the solvers that are used in the experiments. Most of the solvers use the SRR heuristic to create integer solutions and differ from one another by the formulation used to compute the linear relaxation. Moreover, the solvers can be divided into two categories. The first category is made of methods that consider only one decision time step at a time and are given the path that was used by each commodity in the previous time step. These methods can be used to solve horizons with several time steps in a rolling horizon fashion; the paths chosen for each commodity at a time step are given in entry of the method for the next time step. The second category of methods is composed of the multi-time step methods which consider the whole given horizon (10 time steps in our experiments). We start by describing the \textbf{one time step methods}.

\textbf{SRR-arc-node.} In this method, the SRR heuristic is applied to the aggregated arc-node formulation. Commodities with the largest demand are rounded first.

\textbf{SRR-arc-path.} This method consists in applying the SRR heuristic to the extended arc-path formulation where the horizon has been restricted to only one time step. Commodities with the largest demand are rounded first.

\textbf{SRR-restricted.} Similar to "SRR-arc-path"; however, the set of usable paths per commodity is restricted as in Section \ref{sec:restriction}.

\textbf{B\&B-restricted-short/long.} In this method, a standard Branch and Bound procedure is applied to the extended arc-path MILP formulation where the horizon has been restricted to only one time step and the set of usable paths per commodity is restricted. The Branch and Bound procedure is continued until a time limit is reached. In practice, the time limit is set to $\tau = \beta 1.7^{\sqrt{|V|}}$ to obtain a computing time growth similar to the other methods. The $\beta$ coefficient is set to $1/50$ for the short version, $1/2$ for the long version.

We now describe the \textbf{multi time step methods}.

\textbf{SRR-path-sequence.} SRR is applied to the path-sequence formulation. The formulation is solved using the column generation and the pricing scheme presented in the paragraph "Pricing all-in-one" of Section \ref{sec:new_pricing}. During a rounding step, for a commodity, a path-sequence is chosen among the ones present in the linear relaxation through randomized rounding.

\textbf{SRR-path-sequence-restricted.} Similar to "SRR-path-sequence"; however, the set of usable paths per commodity is restricted as in Section \ref{sec:restriction}.

\section{Experimental study}
\label{sec:computational_results}

In this section, we report an experimental comparison of different solvers based on the concepts presented in the previous sections. Several key aspects of the solvers explaining their performance are highlighted through an ablation study. We first introduce how the instances are generated. Then the different solvers and their parameter settings are detailed. Finally, experimental results are presented and discussed. The datasets and the code used in the experimental section of this work are accessible at \url{https://github.com/SuReLI/Dynamic_mcnf_paper_code}. All the code for this work was written in Python 3 and used the commercial solver \citet{gurobi} in its version 8.11. The experiments were made on a server with 48 CPU Intel Xeon E5-2670 2.30GHz, 60 Gbit of RAM, and CentOS Linux 7. The parameter settings for each algorithm are given in \ref{sec:settings}.

\subsection{Instance generation}

In our experiments, we created instances of the dynamic unsplittable flow problem by adapting a method presented in \citep{lamothe2021randomized} for the static problem. All the details of the instance generation process are given in \ref{sec:instance_generation} but we summarize it in the following. The creation of an instance consists of three steps:
\begin{enumerate}
    \item Creation of an initial graph;
    \item Creation of an initial commodity list and initial path for each commodity;
    \item For each time step, modify the graph and the commodity list of the previous time step to create new ones for this time step.
\end{enumerate}
Unless mentioned otherwise in a specific dataset, the capacities of the arcs are set to $10 000$ and the size of the largest commodity possible is set to $1500$. In every dataset, the amount of overflow $B$ allowed in each time step is equal to $1 \%$ of the total demand of the commodities. Moreover, the price of the penalties is set to $1$ which in practice makes the penalties a secondary objective compared to not exceeding the amount of allowed overflow $B$. At each time step, approximately $3\%$ of the commodities change their origin, and each origin of commodities changes $3\%$ of its outgoing arcs. All the instances created contain 10 time steps after the initial one. 

\textbf{The datasets.} Three of the datasets consider graphs of different sizes while the last one considers graphs of fixed sizes but a varying number of commodities. In each dataset, one parameter varies and ten instances are generated for each value of that parameter.
\begin{itemize}
    \item \textit{Grid easy} dataset: This dataset considers grid graphs from 12 nodes to 156 nodes. Many small commodities are created, which makes the instance easier to solve. The capacities of the arcs in the grid are set to 15000 while the capacities of the extra arcs are set to 10000. This also makes the instances easier.
    \item \textit{Grid hard} dataset: This dataset considers grid graphs from 6 nodes to 90 nodes without the adjustments of the previous dataset.
    \item \textit{Random connected} dataset: This dataset considers strongly connected random graphs from 12 nodes to 182 nodes.
    \item \textit{Commodity size} dataset: This dataset considers grid graphs with 42 nodes and uniform arc capacities $c_{et}$ ranging from 1 to 1000 depending on the instance. Moreover, the size of the largest commodity possible is set to $\sqrt{c_{et}}$. This induces a varying number of commodities together with commodities of different sizes compared to the arc capacities.
    \item \textit{Period scaling} dataset: This dataset considers strongly connected random graphs 100 nodes. The number of periods is ranging from 10 to 100.
\end{itemize}

\subsection{Empirical results}
\label{resultat_dynamic}
In this section, each figure presents a subset of algorithms on one of the four datasets for one of three metrics presented below: computing time, overflow ratio, path-change ratio.
In each figure, results from instances created using the same parameters are aggregated. The plotted curves represent the average results on the aggregated instances, while $95 \%$ confidence intervals for the mean are represented as semi-transparent boxes around the main curve. These confidence intervals are created using the statistical method Bootstrap \citep{efron1992bootstrap} with a number of re-samplings equal to 1000.

The path-changes ratio is the ratio of the number of path changes in the solution over the minimum number of path changes of any valid solution. Indeed, even if the capacities of the graph were infinite, the minimum number of path changes in a valid solution is almost always larger than zero. To interpret correctly the solutions yielded by the tested algorithms, we compute the minimum number of path changes achievable when the overflow is not penalized. Note that the value presented in the figures is the path-changes ratio minus one. Indeed, this makes it possible to better highlight (using a logarithmic scale) the quality of the solvers which return solutions whose path-changes ratio is close to one.

The overflow ratio is computed as follows. For each instance, an amount $B$ of overflow is allowed at each time step without penalization. A solution might exceed this allowed amount and have some penalized overflow $o_t$ on certain time steps. The overflow ratio is the total penalized overflow of a solution over the total allowed overflow of an instance: $\frac{\sum_{t \in T} o_t}{|T|B}$. The overflow ratio spans several orders of magnitude while also taking null values. To appropriately display this metric, we use a symmetric logarithmic scale. This scale is logarithmic except around zero where it is linear. 

As for the computing time, it is given in seconds. A time limit of 3 hours has been given to each algorithm. When it is exceeded the corresponding instances are not considered in the average. If only two or fewer instances from a group of instances finished within the time limit, the results of this group are not displayed as we consider the statistical power of the results to be too weak. The points present in the figures without all the instances are:
\begin{itemize}
    \item \textit{Grid hard} dataset: "SRR path-combination" finished only 7 instances with 72 nodes;
    \item \textit{Commodity size} dataset: the first shown points for "SRR path-combination" contains 5 instances while the second only 8;
    \item \textit{Period scaling} dataset: "SRR path-combination" finished 9, 7, and 4 instances with 20, 30, and 40 periods respectively; "SRR path-combination restricted" finished 7 and 4 instances with 40 and 50 periods respectively.
\end{itemize}

\begin{figure}[!ht]
     \centering
     \begin{subfigure}[b]{0.49 \columnwidth}
         \centering
         \includegraphics[width=0.9\textwidth]{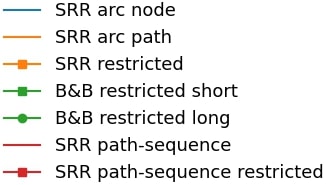}
         \caption{Legend}

     \end{subfigure}
     \hfill
     \begin{subfigure}[b]{0.49 \columnwidth}
         \centering
         \includegraphics[width=\textwidth]{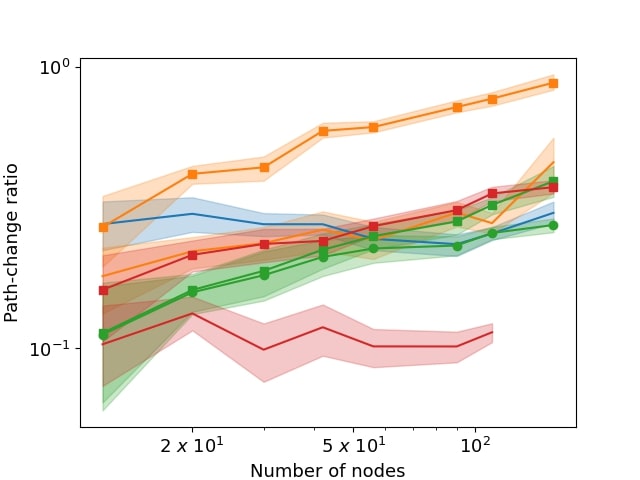}
         \caption{Path-change ratio}
         \label{Important_performance_easy}
     \end{subfigure}
     \begin{subfigure}[b]{0.49 \columnwidth}
         \centering
         \includegraphics[width=\textwidth]{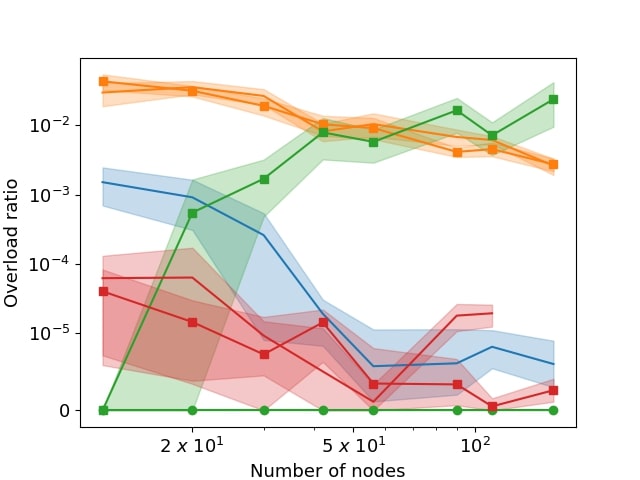}
         \caption{Overload ratio}
         \label{Important_overload_easy}
     \end{subfigure}
     \hfill
     \begin{subfigure}[b]{0.49 \columnwidth}
         \centering
         \includegraphics[width=\textwidth]{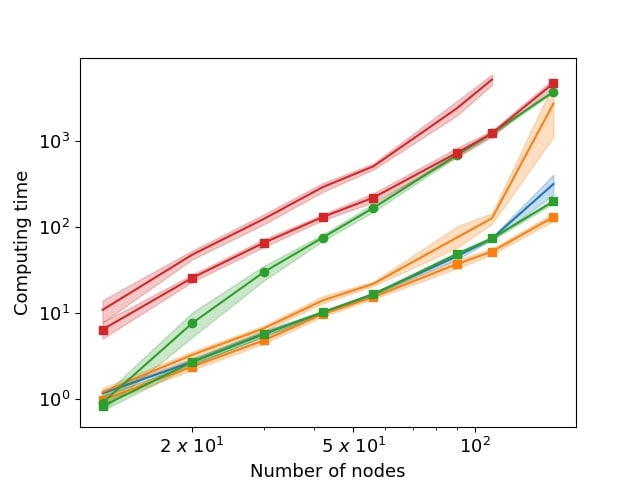}
         \caption{Computing time (s)}
         \label{Important_time_easy}
     \end{subfigure}
        \caption{Results for the \textit{grid easy} dataset}
        \label{Important_easy}
\end{figure}

\begin{figure}[!ht]
     \centering
     \begin{subfigure}[b]{0.49 \columnwidth}
         \centering
         \includegraphics[width=0.9\textwidth]{Legende_en.jpg}
         \caption{Legend}

     \end{subfigure}
     \hfill
     \begin{subfigure}[b]{0.49 \columnwidth}
         \centering
         \includegraphics[width=\textwidth]{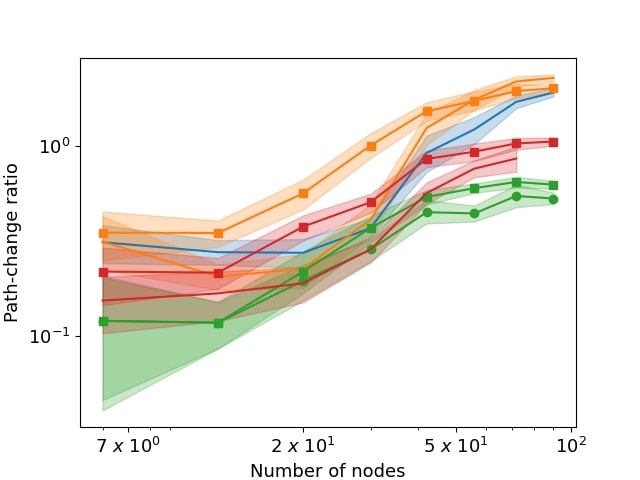}
         \caption{Path-change ratio}
         \label{Important_performance_hard}
     \end{subfigure}
     \begin{subfigure}[b]{0.49 \columnwidth}
         \centering
         \includegraphics[width=\textwidth]{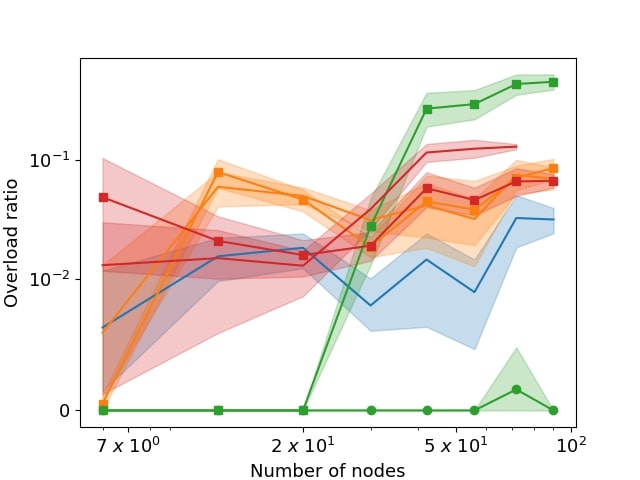}
         \caption{Overload ratio}
         \label{Important_overload_hard}
     \end{subfigure}
     \hfill
     \begin{subfigure}[b]{0.49 \columnwidth}
         \centering
         \includegraphics[width=\textwidth]{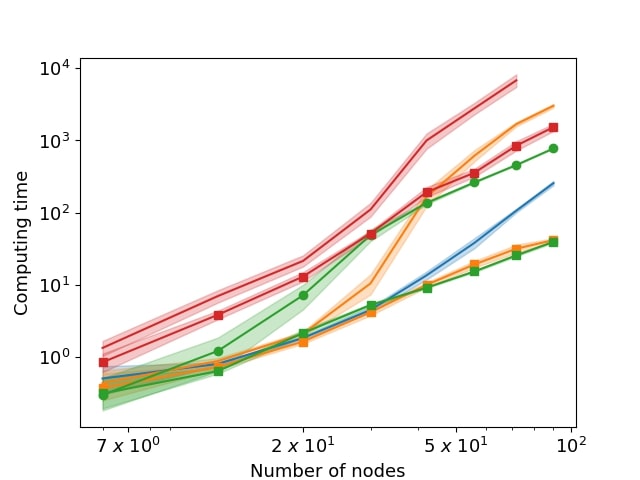}
         \caption{Computing time (s)}
         \label{Important_time_hard}
     \end{subfigure}
        \caption{Results for the \textit{grid hard} dataset}
        \label{Important_hard}
\end{figure}

\textbf{Path set restriction.} In order to analyze the influence of the restriction of the set of allowed paths presented in Section \ref{sec:restriction}, we compare SRR-arc-path to SRR-restricted and SRR-path-sequence to SRR-path-sequence-restricted.  Overall, the path restriction seems to introduce a compromise between computing time and number of path changes. Indeed, in most cases, the computing time is reduced, especially for large instances (typically by 2 to 5 times but sometimes by orders of magnitude on hard instances) while the path-change ratio is increased (typically by 2 to 4 times). This was expected as the path set restriction reduces the search space of the problem which makes the algorithm quicker at the price of not being able to find the optimal solution. As for the overflow, most of the time the difference is not statistically significant.

\begin{figure}[!ht]
     \centering
     \begin{subfigure}[b]{0.49 \columnwidth}
         \centering
         \includegraphics[width=0.9\textwidth]{Legende_en.jpg}
         \caption{Legend}

     \end{subfigure}
     \hfill
     \begin{subfigure}[b]{0.49 \columnwidth}
         \centering
         \includegraphics[width=\textwidth]{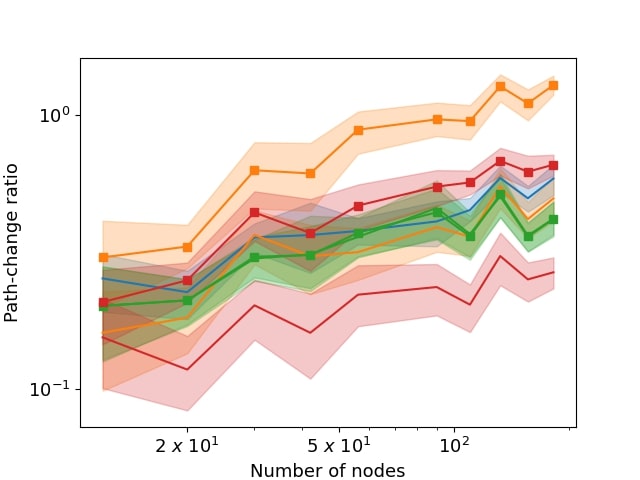}
         \caption{Path-change ratio}
         \label{Important_performance_random}
     \end{subfigure}
     \begin{subfigure}[b]{0.49 \columnwidth}
         \centering
         \includegraphics[width=\textwidth]{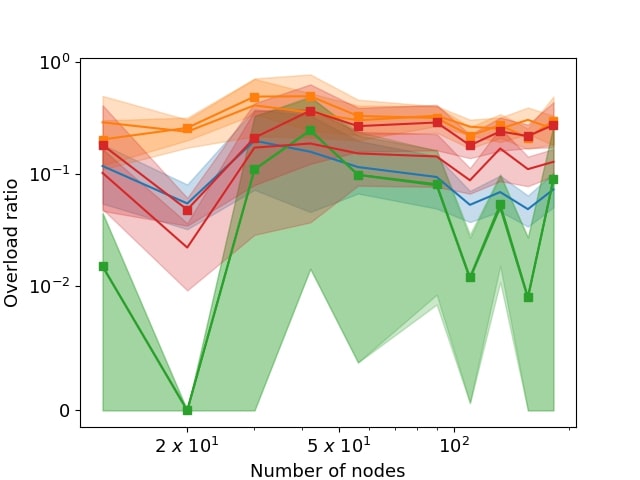}
         \caption{Overload ratio}
         \label{Important_overload_random}
     \end{subfigure}
     \hfill
     \begin{subfigure}[b]{0.49 \columnwidth}
         \centering
         \includegraphics[width=\textwidth]{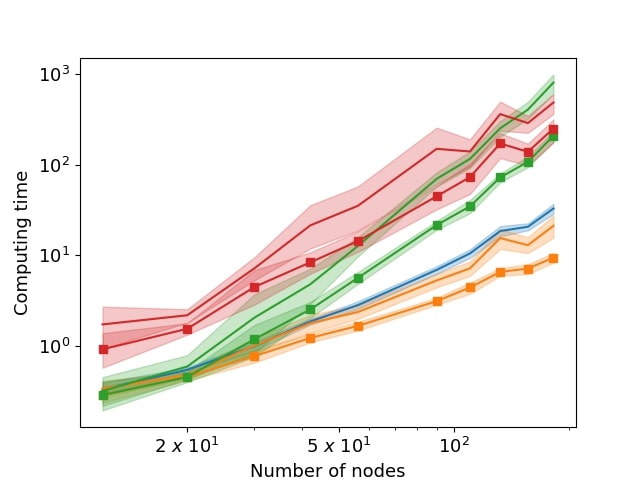}
         \caption{Computing time (s)}
         \label{Important_time_random}
     \end{subfigure}
        \caption{Results for the \textit{random connected} dataset}
        \label{Important_random}
\end{figure}

\begin{figure}[!ht]
     \centering
     \begin{subfigure}[b]{0.49 \columnwidth}
         \centering
         \includegraphics[width=0.9\textwidth]{Legende_en.jpg}
         \caption{Legend}

     \end{subfigure}
     \hfill
     \begin{subfigure}[b]{0.49 \columnwidth}
         \centering
         \includegraphics[width=\textwidth]{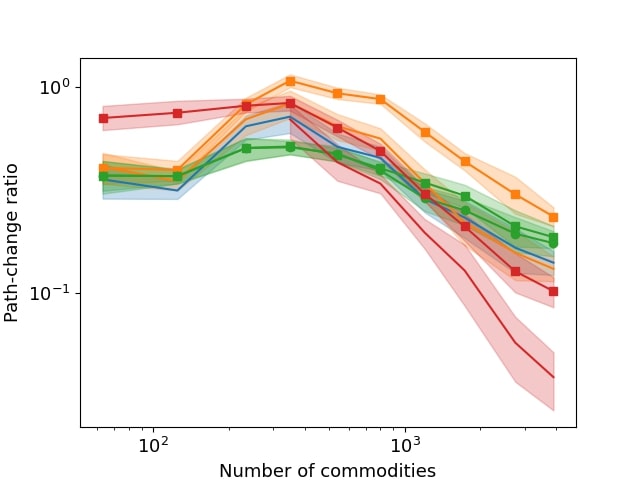}
         \caption{Path-change ratio}
         \label{Important_performance_commodities}
     \end{subfigure}
     \begin{subfigure}[b]{0.49 \columnwidth}
         \centering
         \includegraphics[width=\textwidth]{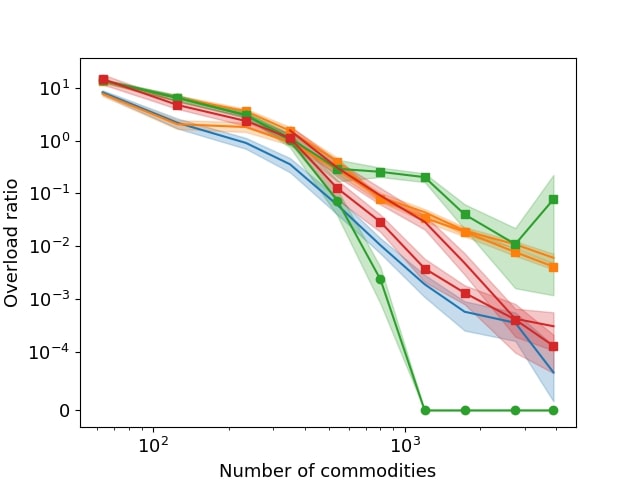}
         \caption{Overload ratio}
         \label{Important_overload_commodities}
     \end{subfigure}
     \hfill
     \begin{subfigure}[b]{0.49 \columnwidth}
         \centering
         \includegraphics[width=\textwidth]{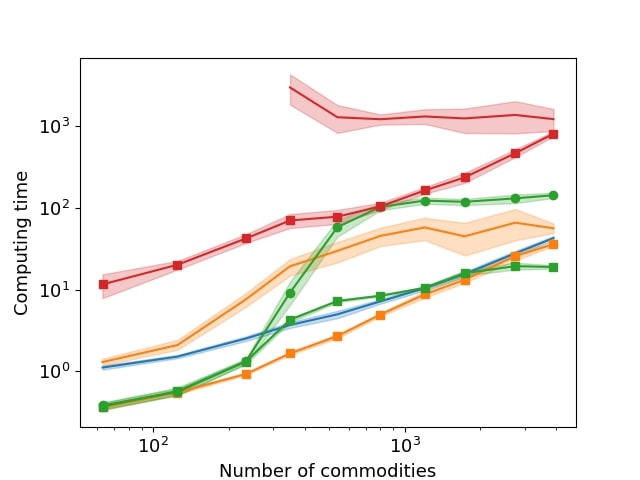}
         \caption{Computing time (s)}
         \label{Important_time_commodities}
     \end{subfigure}
        \caption{Results for the \textit{commodity size} dataset}
        \label{Important_commodity}
\end{figure}

\begin{figure}[!ht]
     \centering
     \begin{subfigure}[b]{0.49 \columnwidth}
         \centering
         \includegraphics[width=0.9\textwidth]{Legende_en.jpg}
         \caption{Legend}

     \end{subfigure}
     \hfill
     \begin{subfigure}[b]{0.49 \columnwidth}
         \centering
         \includegraphics[width=\textwidth]{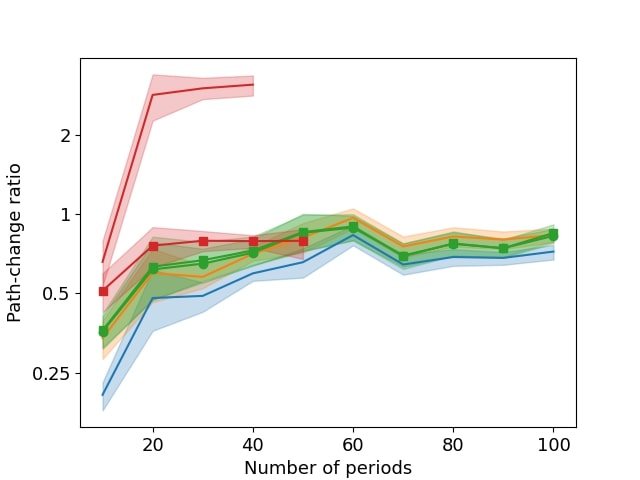}
         \caption{Path-change ratio}
         \label{Important_performance_period}
     \end{subfigure}
     \begin{subfigure}[b]{0.49 \columnwidth}
         \centering
         \includegraphics[width=\textwidth]{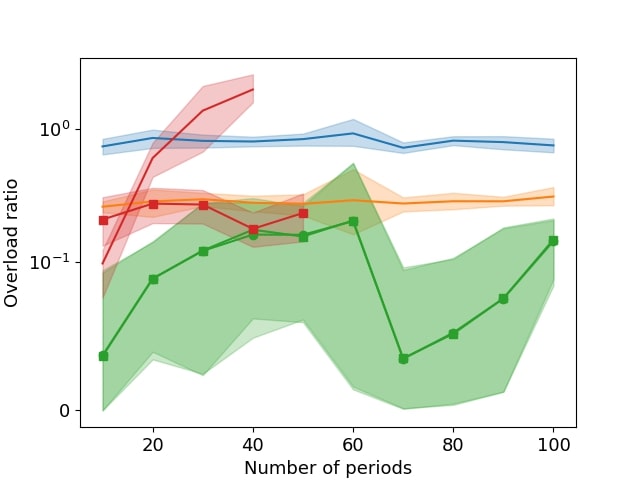}
         \caption{Overload ratio}
         \label{Important_overload_period}
     \end{subfigure}
     \hfill
     \begin{subfigure}[b]{0.49 \columnwidth}
         \centering
         \includegraphics[width=\textwidth]{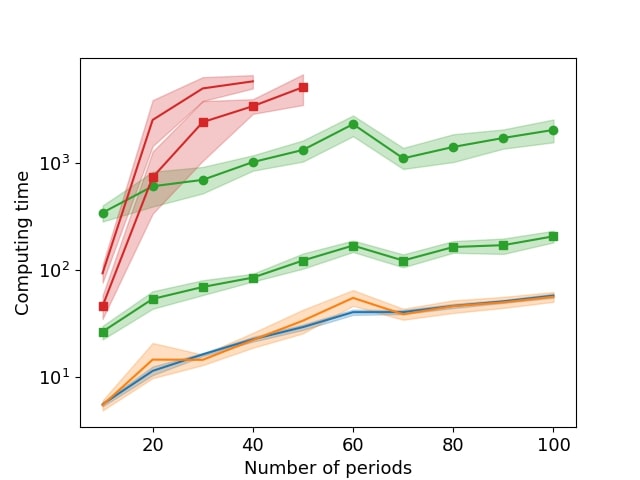}
         \caption{Computing time (s)}
         \label{Important_time_period}
     \end{subfigure}
        \caption{Results for the \textit{period scaling} dataset}
        \label{Important_period}
\end{figure}

\textbf{Aggregated arc-node versus arc-path formulation.} Two formulations were used to compute exact linear relaxations in methods considering one time step at a time: the arc-path formulation in SRR-arc-path and the aggregated arc-node formulation in SRR-arc-node. Although both methods yield solutions with similar path-change ratio SRR-arc-node runs, most of the time, faster than SRR-arc-path and yield a lower overflow ratio when the difference is statistically significant. Note that both methods exploit the fact that, in our instances, commodities originate from a small number of sources. SRR-arc-node uses aggregate variables while SRR-arc-path makes fewer calls to Dijkstra's shortest path algorithm when generating the paths used by each commodity. Overall, when grouping the commodities is possible, using an aggregated arc-node model yields better results than using an arc-path model. However, SRR-arc-node may become intractable when the number of sources is large because of a large number of variables.

\textbf{Scaling to longer horizons.} According to Figure \ref{Important_period}, the one time step methods scale rather well into longer horizons. Indeed, their computing time appears to grow linearly with the number of periods and the quality of the returned solutions seems to be quite stable. On the other hand, the computing time of the multi-timestep methods grows very quickly. Its growth is even larger than what appears at first glance in Figure \ref{Important_time_period} because the last points display only the results for the instances that were solved in the 3 hours time limit. However, for the instances that were solved, the objective value of the solutions appears quite stable. This can be seen for "SRR path-combination restricted" directly on the figures however this is probably also the case for "SRR path-combination". We hypothesize that the rapid degradation of solution quality for this method is due to the fact that the computation of a good linear relaxation needs more column generation iterations than we could allow this method to have while still finishing in the 3 hours time limit.

\textbf{General comments.} As the number of nodes increases all methods need more time to compute solutions with a higher path-change ratio. However, the overload ratio is either stagnating or decreasing. As for the number of commodities, similarly to the static unsplittable flow problem, the solution method still requires more computation time but they are able to find solutions of better quality both in terms of overload ratio and path-change ratio. Although the methods considering several time steps are significantly slower than their one-time step counterparts, most of the time they deliver solutions of better quality. Finally, applying a Branch and Bound method to the arc-path formulation with restricted usable paths appears to be one of the best compromises between computing time and solution quality among the tested algorithms. Branch and Bound methods can also benefit from the presence of several CPUs even though such results are not presented here to ensure a fair comparison with the other methods. The fact that an initial path is given and preferable to the other paths seems to enable good results from the internal heuristics of the commercial MILP solvers. Indeed, such performance could not be reproduced on the static unsplittable flow problem where no path is preferred among the set of valid paths.

\section{Conclusion}
\label{sec:conclusion}
In this paper, several new methods to solve medium to large instances of the dynamic unsplittable flow problem were presented. In particular, new formulations were introduced which model either the integer problem or its linear relaxation. Moreover, we introduced new approaches to solve the pricing problem for the formulation introduced by \citet{gamvros2012multi}. These methods do not rely on k-shortest path computations and achieve a polynomial worst-case complexity.

The formulations were embedded in matheuristic solvers for the dynamic unsplittable flow problem which were compared on several benchmarks of instances. Several key aspects of the solvers explaining their performances were also highlighted. Overall, considering only one time step at a time and a small set of allowed paths for each commodity gives a good trade-off between solution quality and computing time when the formulation is solved using a commercial solver. However, the methods that yield the solutions with the best quality revolve around the path-sequence formulation which considers several time steps at a time.

As for future works, although meta-heuristics were not considered in this work and were not yet used in the literature for the dynamic unsplittable flow problem, they could be investigated as an alternative to methods based on (mixed-integer) linear formulations.

\bibliographystyle{plainnat}
\bibliography{biblio.bib}

\newpage

\appendix

\section{Computational results for the extended arc-node formulation}
\label{appendix_arc_node}

In this section, we present the extended arc-node formulation and give computational results on very small instances.
\begin{itemize}
    \item $x_{et}^k$ decides if commodity $k$ uses arc $e$ at time step $t$ to send its flow;
    \item $o_{et}$ represents the overflow on arc $e$ at time step $t$;
    \item $o_t$ represents the amount of overflow that exceeds, at time step $t$, the amount $B$ of unpenalized overflow;
    \item $n_{t}^k$ takes value one if commodity $k$ makes a path-change between time step $t$ and time step $t-1$.
\end{itemize}

\begin{subequations}
\begin{alignat}{2}
&\min_{x_{et}^k,n_{t}^k, o_{et}, o_t} ~ \alpha \sum_{k \in K, t \in T} n_{t}^k + \sum_{t \in T} o_t & \qquad &\\
&\text{subject to} &      & \notag \\
&\sum_{e \in E^+(v)} \hspace{-0.3cm} x_{et}^k ~ \hspace{-0.2cm} - \hspace{-0.3cm} \sum_{e \in E^-(v)} \hspace{-0.3cm} x_{et}^k = \delta^{O^k}_v - \delta^{D^k}_v & & \forall k' \in K', ~ \forall v \in V \label{eq: flow_const}\\
&\sum_{k \in K} x_{et}^k ~ d^k \leq c_{et} + o_{et} & & \forall e \in E_t, ~ \forall t \in T \label{eq: capacity_const_2}\\
&\sum_{e \in E_t} o_{et} \leq B + o_t & &  \forall t \in T \label{eq: capacity_const_tot_2}\\
& x_{et}^k - x_{e,t-1}^k \leq n_{t}^k & & \forall e \in E, ~ \forall k \in K, ~ \forall t \in T \setminus \{0\} \label{eq: delta_2}\\
& x_{et}^k \in \{0,1\}, ~ n_{t}^k \in \mathbb{R}^+  & & \forall e \in E, ~ \forall k \in K, ~ \forall t \in T \label{eq: integrity_2}\\
& o_{et} \in \mathbb{R}^+, ~ o_t \in \mathbb{R}^+ & & \forall e \in E_t, ~ \forall t \in T
\end{alignat}
\end{subequations}
The objective function is composed of two terms: the sum of the path-change penalties and the sum over the time steps of the overflow exceeding the threshold $B$. The equations \ref{eq: flow_const} are the flow conservation constraint ensuring that the arcs selected by the variables $x_{et}^k$ represent paths. The equations (\ref{eq: capacity_const_2}) and (\ref{eq: capacity_const_tot_2}) are the capacity constraints: they ensure that $o_{et}$ represents the overflow on arc $e$ at time step $t$ and that $o_t$ represents the amount of overflow that exceeds the threshold $B$ at time step $t$. The equations (\ref{eq: delta_2}) ensure that $n_{t}^k$ takes value 1 when a change of path occurs for commodity $k$ between time step $t$ and time step $t-1$. The variables $x_{et}^k$ being binary in the equations \eqref{eq: integrity_2} ensures that the flow is unsplittable.

As it can be seen in Table \ref{results_arc_node}, an exact solver based on an arc-node formulation does not scale well and struggles even for the smallest instances of the \textit{grid hard} dataset. It is to be noted that for graphs with 12 nodes, 6 out of 10 instances were not solved in the time limit of 3 hours, and the remaining four instances were solved in an average of 1h30. One can also note that the arc-node formulation finds solutions with a zero overload ratio which comes at the cost of a lot more path changes than the other methods. These are indeed better solutions because the overflow after threshold B is very highly penalized in our instances. 

\begin{table}
    \centering
    \begin{tabular}[c]{|l|c|c|c|c|c|c|}
    \hline
     & \multicolumn{2}{c|}{Path-change ratio}   & \multicolumn{2}{c|}{Overload ratio}   & \multicolumn{2}{c|}{Computing time (s)}  \\ \hline
    Instance size (nodes) & 6 & 12 & 6 & 12 & 6 & 12\\ \hline
    Arc-node formulation & 1.95   & 2.85 & 0.0 & 0.0 & 557 & 5467 \\
    SRR arc node & 0.15 & 0.28 & 0.32 & 0.16 & 0.33 & 1.15 \\
    SRR path combination & 0.23 & 0.39   & 0.013 & 0.016 & 1.32 & 14.3 \\
    \hline
    \end{tabular}
    \caption{Computational results for the arc-node formulation on small instances : averaged over 10 instances for each size.}
    \label{results_arc_node}
\end{table}

\section{Proof of the equivalence of the strength of the linear relaxations of the extended arc-path formulation and path-sequence formulation}
\label{appendix_proof}

We recall here the proposition stating the equivalence of strength of the linear relaxation of the two formulations and give a formal proof.

\begin{proposition}
Let $V_1$ be the value of the linear relaxation of the extended arc-path formulation and $V_2$ be the value of the linear relaxation of the path-sequence formulation, then $V_1 = V_2$.
\end{proposition}

\begin{proof}

Let $R$ be the polyhedron induced by the two capacity constraints \eqref{eq: capacity_const} and \eqref{eq: capacity_const_tot} and let $Q_1$ be the polyhedron:
$$Q_1 = \{x_{pt}^k \in [0, 1], ~ n_{pt}^k \in \mathbb{R}^ +, \text{ satisfying constraints } \eqref{eq: convex} \text{ and } \eqref{eq: delta_1}\}.$$
The solution space of the linear relaxation of the extended path formulation is $ R \cap Q_1 $. Meanwhile, the path-sequence formulation is obtained by applying a Dantzig-Wolfe decomposition to the polyhedron $Q_1$. Thus, the solution space of its linear relaxation is $ R \cap Q_2 $ where $Q_2$ is the convex envelope of integer points of $Q_1$. In order to show that the value of both linear relaxations is the same, we will show that $ Q_1 = Q_2 $.

The inclusion $ Q_2 \subseteq Q_1 $ is implied by the fact $Q_2$ is the result of the Dantzig-Wolfe decomposition of $Q_1$. Thus, let us now show $ Q_1 \subseteq Q_2 $.

To that end, consider an assignment of the variables $ x_{pt}^k $ and $n_{pt}^k$ satisfying constraints \eqref{eq: convex} and \eqref{eq: delta_1} with the variables $n_{pt}^k$ as small as possible. Note that the variable $n_{pt}^k$ will take the value $\max(0, x_{pt}^k - x_{p,t-1}^k)$. We will show there exists an assignment of the variables $x_s^k$ of the path-sequence formulation inducing the same flow distribution and the same number of path-change penalties. We now present an algorithm to compute such an assignment. This construction is carried out separately for each commodity and, therefore, we consider only one commodity in the following.

In what follows, let us call $ \Pi_{pt}^k $ the set of path-sequences constructed by the algorithm for commodity $ k $ and using the path $ p $ at the time step $ t $.  Moreover, let us denote $ \lambda(s) $ the coefficient of a path-sequence $ s $ in a convex combination of path-sequences and let $ \lambda (\Pi) = \sum_{s \in \Pi} \lambda(s) $.

\textit{Algorithm outline.} The algorithm constructs a set of path-sequences as well as their coefficients in the convex combination. First, all path-sequences are initialized with a path for the first time step, then they are all extended with a path for the second time step, and so on. This construction has two underlying objectives. First, we want $\lambda(\Pi_{pt}^k) = x_{pt}^k $, so that path-sequences represent the flow distribution induced by $ x_{pt}^k $. Second, we extend, if possible, the path-sequences with the path they already end with so that a minimum number of path-changes is induced.

We give an example of execution of the above algorithm in Table \ref{example_proof}. In this example, the distribution of the flow represented by the variables $ x_{pt}^k$ is given in Table \ref{example_x_values}. In the other tables, we present the execution of the algorithm. In these tables, the column $ t $ contains the path-sequences - indexed by the letter $ s $ - created by the algorithm during the extensions of time step $ t $.

\begin{table}
    \begin{subtable}[h]{0.45\textwidth}
        \centering
        \begin{tabular}[c]{|c|ccc|}
        \hline
        \backslashbox{p}{t} & 1   & 2   & 3  \\ \hline
        a & 0.5 & 0.3 & 0.4 \\
        b & 0.2 & 0   & 0 \\
        c & 0   & 0.5 & 0.3 \\
        d & 0.3 & 0.2 & 0.3\\
        \hline
        \end{tabular}
        \caption{Flow distribution: value of the $x_{pt}^k$ variables for one commodity}
        \label{example_x_values}
    \end{subtable}
    \hfill
    \begin{subtable}[t]{0.45\textwidth}
        \centering
        \begin{tabular}{|c|c|c|c|}
        \hline
        \backslashbox{$s$}{$t$} & 1 & 2 & 3 \\ \hline
        1 & $a$ | 0.5 & ~~~~~~~~ & ~~~~~~~~~~ \\ \hline
        2 & $b$ | 0.2 &  &  \\ \hline
        3 & $d$ | 0.3 &  &  \\ \hline
        4 &           &  &  \\ \hline
        5 &           &  &  \\ \hline
        6 &           &  &  \\ \hline
        7 &           &  &  \\ \hline
        8 &           &  &  \\ \hline
        \end{tabular}
        \caption{Initialization: a sequence containing a single path is created for each path used at the first time step}
        \label{example_lambda_values1}
     \end{subtable}

    \begin{subtable}[t]{0.45\textwidth}
        \centering
        \begin{tabular}{|c|c|c|c|}
        \hline
        \backslashbox{$s$}{$t$} & 1 & 2 & 3 \\ \hline
        1 & $a$ | 0.5 & $aa$ | 0.3 & ~~~~~~~~~~ \\ \hline
        2 & $b$ | 0.2 &            &  \\ \hline
        3 & $d$ | 0.3 & $dd$ | 0.2 &  \\ \hline
        4 &           &            &  \\ \hline
        5 &           &              &  \\ \hline
        6 &           &            &  \\ \hline
        7 &           &            &  \\ \hline
        8 &           &            &  \\ \hline
        \end{tabular}
        \caption{Prioritized extensions for $t = 2$: the sequences are extended with the path they already end with. The coefficients are decreased so that the sum of the coefficients of the sequences ending with a path $ p $ does not exceed $x_{p2}^k$}
        \label{example_lambda_values2}
    \end{subtable}
    \hfill
    \begin{subtable}[t]{0.45\textwidth}
        \centering
        \begin{tabular}{|c|c|c|c|}
        \hline
        \backslashbox{$s$}{$t$} & 1 & 2 & 3 \\ \hline
        1 & $a$ | 0.5 & $aa$ | 0.3 & ~~~~~~~~~~ \\ \hline
        2 & $b$ | 0.2 &            &  \\ \hline
        3 & $d$ | 0.3 & $dd$ | 0.2 &  \\ \hline
        4 &           & $ac$ | 0.2 &  \\ \hline
        5 &           & $bc$ | 0.2 &  \\ \hline
        6 &           & $dc$ | 0.1 &  \\ \hline
        7 &           &            &  \\ \hline
        8 &           &            &  \\ \hline
        \end{tabular}
        \caption{Other extensions for $t = 2$: the sequences whose coefficient has been reduced (here all of them) are extended with paths where there is still some flow to allocate (here only the path $ c $)}
        \label{example_lambda_values3}
     \end{subtable}
     
    \begin{subtable}[t]{0.45\textwidth}
        \centering
        \begin{tabular}{|c|c|c|c|}
        \hline
        \backslashbox{$s$}{$t$} & 1 & 2 & 3 \\ \hline
        1 & $a$ | 0.5 & $aa$ | 0.3 & $aaa$ | 0.3 \\ \hline
        2 & $b$ | 0.2 &            &             \\ \hline
        3 & $d$ | 0.3 & $dd$ | 0.2 & $ddd$ | 0.2 \\ \hline
        4 &           & $ac$ | 0.2 & $acc$ | 0.2 \\ \hline
        5 &           & $bc$ | 0.2 & $bcc$ | 0.1 \\ \hline
        6 &           & $dc$ | 0.1 &             \\ \hline
        7 &           &            &             \\ \hline
        8 &           &            &             \\ \hline
        \end{tabular}
        \caption{Prioritized extensions for $t = 3$: the sequences are extended with the path they already end with. Note in particular that the coefficients of $ bcc $ and $ dcc $ are decreased (down to 0 for $ dcc $) so that $\lambda(acc) + \lambda(bcc) + \lambda(dcc) = x_{c3}^k$}
        \label{example_lambda_values4}
    \end{subtable}
    \hfill
    \begin{subtable}[t]{0.45\textwidth}
        \centering
        \begin{tabular}{|c|c|c|c|}
        \hline
        \backslashbox{$s$}{$t$} & 1 & 2 & 3 \\ \hline
        1 & $a$ | 0.5 & $aa$ | 0.3 & $aaa$ | 0.3 \\ \hline
        2 & $b$ | 0.2 &            &             \\ \hline
        3 & $d$ | 0.3 & $dd$ | 0.2 & $ddd$ | 0.2 \\ \hline
        4 &           & $ac$ | 0.2 & $acc$ | 0.2 \\ \hline
        5 &           & $bc$ | 0.2 & $bcc$ | 0.1 \\ \hline
        6 &           & $dc$ | 0.1 &             \\ \hline
        7 &           &            & $bca$ | 0.1 \\ \hline
        8 &           &            & $dcd$ | 0.1 \\ \hline
        \end{tabular}
        \caption{Other extensions for $t = 3$: the sequences whose coefficient has been reduced (here $ bc $ and $ dc $) are extended with paths where there is still some flow to be allocated (here the paths $ a $ and $ d $)}
        \label{example_lambda_values5}
     \end{subtable}
     \caption{Example of execution of the algorithm used in the proof}
     \label{example_proof}
\end{table}

\textit{Initialization of the convex combination of the path-sequences.} For each path $ p \in P^k_1$, the algorithm creates a path-sequence containing only path $ p $ with a coefficient equal to $ x_{ p1}^k$ and adds it to $ \Pi_{p1}^k $.

Then, for each following time step $ t $, the algorithm proceeds as follows.

\textit{Prioritized extensions.} The algorithm first begins by extending as much as possible the path-sequences with the path they already end with. Let us consider a path $ p \in P^k_t$ and distinguish two cases. 
\begin{itemize}
    \item If $x_{pt}^k \geq x_{p,t-1}^k$ then all the path-sequences in $\Pi_{p, t-1}^k $ are extended with path $ p $ and added to $ \Pi_{pt}^k $. In this case, at this point, $\lambda(\Pi_{pt}^k) = x_{p,t-1}^k \leq x_{pt}^k$.
    \item If $x_{pt}^k < x_{p,t-1}^k$, the algorithm finds a subset $ \Pi \subset \Pi_{p,t-1}^k$ and a path-sequence $s \in \Pi_{p,t-1}^k$ such that: $s \notin \Pi$ and $\lambda(\Pi) \leq x_{pt}^k \leq \lambda (\Pi) + \lambda(s)$. Then, the path-sequence $ s $ is "split" into two path-sequences; $ s_1 $ which is given $ x_{pt}^k - \lambda(\Pi) $ as coefficient and $ s_2 $ which is given $\lambda(s) - \lambda(s_1)$ as coefficient. This choice of coefficients will ensure that $\lambda(s_1) + \lambda(s_2) = \lambda(s)$ and $\lambda(\Pi_{pt}^k) = \lambda(\Pi) + \lambda (s_1) = x_{pt}^k$. Indeed, the path-sequences in $ \Pi \cup \{s_1\} $ are extended with path $ p $ and added to $ \Pi_{pt}^k$ while $s_2$ remains in $\Pi_{p,t-1}^k$.
\end{itemize}
 If this extension process is repeated for each path $ p $, then a maximum amount (in the sense of the coefficients) of path-sequences is extended with the path they already end with.

\textit{Other extensions.} After the prioritized extensions, the algorithm considers the path-sequences that have not yet been extended for the time step $ t $. The algorithm extends the first of these sequences $ s $ with any path $ p $ such that $\lambda(\Pi_{pt}^k) < x_{pt}^k$. Once again, if $\lambda(\Pi_{pt}^k) + \lambda(s) > x_{pt}^k$ then to obtain $\lambda(\Pi_{pt}^k) = x_{pt}^k$, the sequence is divided into two parts with the same coefficients as in the prioritized extensions. The first part is extended with the path $ p $ before being added to $\Pi_{pt}^k$. The other part will be extended later in the algorithm with another path and remains in $ \Pi_{p, t-1}^k $. The algorithm then extends another path-sequence. Once all the sequences have been extended, the extensions can begin for the next time step.

In the above algorithm, the path-sequences are extended so that the sum of the coefficients of the path-sequences using the path $ p $ at the time step $ t $ is $ x_{pt}^k $. Thus, the convex combination implies the same flow distribution as the variables $ x_{pt}^k $. It remains to verify that the two representations of the flow distribution imply the same number of path-changes. In the algorithm, at the time step $ t $, the sum of the coefficients of the path-sequences ending with the path $ p $ and extended by the path $ p $ is $\min(x_{pt}^k, x_{p, t-1}^k)$. All other path-sequences involve path-changes whose sum amounts to $x_{pt}^k - \min(x_{pt}^k, x_{p,t-1}^k)$ which is equal to $\max(0, x_{pt}^k - x_{p,t-1}^k)$. This is exactly the value taken by the $ n_{pt}^k $ variables that model the path-changes in the extended arc-path formulation. Thus, the number of path penalties is also the same.

With the above algorithm, any assignment of the variables $ x_{pt}^k $ can be translated into a convex combination of path-sequences inducing the same number of path changes. This means that $ Q_1 \subseteq Q_2 $. Thus, the linear relaxation of the extended path formulation is as strong as the linear relaxation of the path-sequence formulation. Thus, we have finally shown that the two linear relaxations are equivalent.
\end{proof}

\section{Parameter settings}
\label{sec:settings}

In this section, we describe some additional hyper-parameters of our algorithms along with their value in the experiments.

\textbf{SRR threshold $\mathbf{\theta}$} This threshold decides how often the linear relaxation is updated and we searched for a value giving a good trade-off between performance and computation time. Its value is fixed to $\theta = |V|$ for methods based on the path-sequence formulation and $|V|/10$ otherwise.

\textbf{Flow penalization} An small additional cost was added to the variables representing the flow in the models. This helps the algorithm focus on shorter paths. In our preliminary tests, this penalization appeared to diminish the computing time and increase the quality of the returned solution. Thus, all the algorithms tested use this penalization. For path variables representing a path $p$, the additional cost is $|E(p)| \epsilon$ where $|E(p)|$ is the number of arcs in path $p$. For path-sequences variables representing a path-sequence $(p_1, ..., p_T)$, the additional cost is $\sum_{t \leq T} |E(p_t)| \epsilon$. For the aggregated arc-node model, the flow variable associated with each arc is given a cost of $\epsilon$. In all cases, the constant $\epsilon$ was set to the small value of $10^{-4}$ in our experiments.

\textbf{Variable deletion} When generating variables in the arc-path model or the path-sequence model, the models tend to accumulate a large number of unnecessary variables (variables not used in the current optimal solution) which increases unnecessarily the resolution time of the model. To prevent this, each time the model is solved, each variable not used in the basis of the optimal solution is deleted with a fixed probability. The value of this probability is fixed at $0.3$ as it appeared to give the best results in our preliminary tests.

\textbf{Size of the restricted path sets} Some of the solvers consider only a restricted number of paths for each commodity as explained in Section \ref{sec:restriction}. The number $\kappa$ of k-shortest paths computed is set to 4 in our tests.

\textbf{Number of column generation iterations} The solvers based on the path-sequence formulation do not perform all the column generation iterations required to obtain the optimal solution of the linear relaxation. This decreases the overall computing time of the algorithms. Before the first randomized rounding step, 20 column generation iterations are made to obtain a good approximation of the linear relaxation. However, between the randomized rounding steps, only 3 column generation iterations are made. Making more iterations appeared unnecessary in our preliminary tests.

\textbf{Pricing scheme for the path-sequence formulation} To solve the path-sequence formulation, a column generation process must be used which relies on a pricing method able to compute the variables with the most negative reduced cost. Three pricing schemes have been presented in Section \ref{sec:solving_path-sequence_formulation}. The one introduced by \citet{gamvros2012multi} is not applicable in our variation of the problem due to its limitations presented in Section \ref{sec:pricing_gamvros}. Although the other two methods do not have the same theoretical complexity, preliminary results showed that choosing one over the other did not significantly impact the computing time in our instance. The pricing scheme used in the results is the one presented in the paragraph "Pricing all-in-one" of Section \ref{sec:new_pricing} which has the lowest complexity.

\textbf{Branch and Bound settings} The parameter MIPFocus of \citet{gurobi} was set to 1 to encourage the solver to find high-quality solutions quickly. Moreover, to keep a fair comparison with the other solvers, the solvers using the Branch and Bound of \citet{gurobi} have been restricted to using only one CPU. This prevented the Branch and Bound solvers to use the dozens of CPUs available on the server where the tests were made.

\section{Instance generation}
\label{sec:instance_generation}

The creation of an instance of the dynamic unsplittable flow problem consists of four steps:
\begin{enumerate}
    \item Creation of the initial graph;
    \item Creation of the initial commodity list;
    \item Creation of the initial path affected to each commodity;
    \item For each time step, modify the graph and the commodity list of the previous time step to create new ones for this time step.
\end{enumerate}

The initial graph and commodity list of the instances are created with a method for instances of the static unsplittable flow problem presented in \citep{lamothe2021randomized}. We describe this method in the two following paragraphs.

\textbf{Graph types.} We consider two types of graphs: strongly connected random graphs and grid graphs. For strongly connected random graphs, we ensure that the graph has only one strongly connected component and then control the average degree of the graph (it cannot be less than 2). The probability of an arc's existence is fixed to $5/|V|$ and the probability of a node being an origin is $1/10$. For the grid type, a $n$-sized graph is an $n \times n$ toric grid with $n$ additional nodes. These additional nodes are each randomly connected to $2n$ nodes on the grid and serve as origins of the flow. Unless mentioned otherwise, the arc capacities are $10^4$.

\textbf{Commodity list.} The initial commodity list is created as follows. For each commodity:
\begin{itemize}
    \item choose a destination node $d$;
    \item choose an origin $o$ which can access $d$ within the remaining capacities;
    \item compute a random simple path $p$ from $o$ to $d$ using a depth-first search where the visiting order of newly discovered nodes is random; we denote $c_p$ the remaining capacity on path $p$;
    \item choose a demand level $D$; for this choice, the parameter $\hat{D}_{\max}$ defines the maximum possible demand of a commodity, and $U(x)$ is a random integer variable uniformly drawn in $[1, x]$; we use two different methods to choose the demand level because it impacts the difficulty of the instance; either $D = \min(c_p, U(\hat{D}_{\max}))$ or $D = U(\min(c_p, \hat{D}_{\max}))$; Unless stated otherwise $\hat{D}_{\max} = 1500$;
    \item decrease the used capacity on the path $p$ by $D$;
    \item add $(o, d, D)$ to the list of created commodities;
    \item repeat until no commodity can be added without breaking the capacity constraints.
\end{itemize}

Note that a commodity list created this way can always be routed within the arc capacities by using, for each commodity, the path $p$ used to create the commodity. However, this is only true for the first time step. In the subsequent time step, since changes are made to the commodities and the graph, there is no guarantee that all the commodities can still be routed inside the capacities. Note that having commodities that fit in the capacities while fully congesting the graph tends to create instances that are rather hard to solve. Indeed, if there are too few commodities, most solutions have zero overflow and the instance is easy to solve. On the other hand, if there is way too much demand to fit in the capacities, overflow will be created whichever path is chosen for each commodity and, again, the instance is easy to solve.

\textbf{Initial path list.} The list of initial paths is chosen to be the list of the paths $p$ used to create the commodities.

\textbf{Subsequent time steps} An instance of the dynamic unsplittable flow problem is a sequence of unsplittable flow problems representing different time steps. Each time step introduces a few changes to the problem: some commodities change their origin or their destination, and some arcs are added or deleted. The following changes are made between each time step in our tests:
\begin{itemize}
    \item For each commodity $(o, d, D)$, with probability $\mu$, the destination $d$ of the commodity is replaced by a node $u$ that is not the origin of a commodity and such that $(d, u) \in E_t$. If no such node exists, the destination remains unchanged.
    \item For each arc $e = (o, u)$ going out of an origin $o$, with probability $\mu$, arc $e$ is replaced by an arc $e' = (o, v)$ where $v$ is a node that is not the origin of a commodity such that $(v, u) \in E_t$. If no such node exists, no change occur.
\end{itemize}

The probability $\mu$ of making a change is set to $\mu = 0.03$. Note that by changing only the destinations of the commodities, we keep the number of different origins low. Also note that, with these arc changes, a strongly connected graph remains strongly connected. All the instances created contain 10 time steps after the initial one.

\textbf{The datasets.} Four different datasets are used during the experiments. Three of them consider graphs of different sizes while the last one considers graphs of fixed sizes but a varying number of commodities. Unless mentioned otherwise in a specific dataset, the capacities of the arcs are set to $10 000$ and the size of the largest commodity possible is set to $1500$. In every dataset, the amount of overflow $B$ allowed in each time step is equal to $1 \%$ of the total demand of the commodities. Moreover, the price of the penalties is set to one which in practice makes the penalties a secondary objective compared to not exceeding the amount of allowed overflow $B$. In each dataset, one parameter varies and ten instances are generated for each value of that parameter.
\begin{itemize}
    \item \textit{Grid easy} dataset: This dataset considers grid graphs from 12 nodes to 156 nodes. When choosing the demand level of a commodity, the formula $D = U(\min(c_p, \hat{D}_{\max}))$ is used. This tends to create many small commodities, which makes the instance easier to solve. The capacities of the arcs in the grid are set to 15000 while the capacities of the extra arcs are set to 10000. This also makes the instances easier.
    \item \textit{Grid hard} dataset: This dataset considers grid graphs from 6 nodes to 90 nodes. When choosing the demand level of a commodity, the formula $D = \min(c_p, U(\hat{D}_{\max}))$ is used.
    \item \textit{Random connected} dataset: This dataset considers strongly connected random graphs from 12 nodes to 182 nodes. When choosing the demand level of a commodity, the formula $D = \min(c_p, U(\hat{D}_{\max}))$ is used.
    \item \textit{Commodity size} dataset: This dataset considers uniform arc capacities $c_{et}$ ranging from 1 to 1000 depending on the instance, while the parameter $\hat{D}_{\max}$ is set to $\sqrt{c_{et}}$. This induces a varying number of commodities together with commodities of different sizes compared to the arc capacities. All considered graphs are grid graphs with 42 nodes. When choosing the demand level of a commodity, the formula $D = \min(c_p, U(\hat{D}_{\max}))$ is used.
    \item \textit{Period scaling} dataset: This dataset considers strongly connected random graphs 100 nodes. The number of periods is ranging from 10 to 100. When choosing the demand level of a commodity, the formula $D = \min(c_p, U(\hat{D}_{\max}))$ is used.
\end{itemize}

\begin{figure}
    \begin{subfigure}[b]{0.49 \columnwidth}
        \centering
        \includegraphics[width=\textwidth]{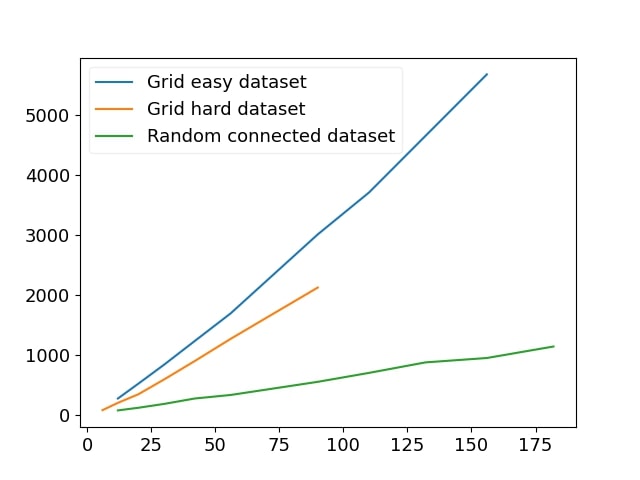}
        \caption{Average number of commodities versus number of nodes}
    \end{subfigure}
    \begin{subfigure}[b]{0.49 \columnwidth}
        \centering
        \includegraphics[width=\textwidth]{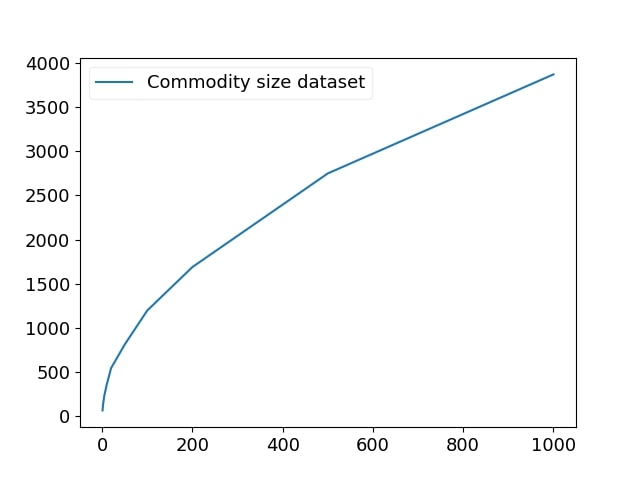}
        \caption{Average number of commodities versus arc capacities}
    \end{subfigure}
    \caption{Average number of commodities for each dataset}
    \label{commodity_figures}
\end{figure}

The average number of commodities relative to the changing parameter of each dataset is given in Figure \ref{commodity_figures}.

\end{document}